\documentclass[a4paper]{article}
\usepackage[margin=30mm]{geometry}
\usepackage{amsmath}
\usepackage{amsfonts}
\usepackage{amssymb}
\usepackage{graphicx}
\usepackage{verbatim}
\def\tr{{\rm{tr}}}
\newtheorem{theorem}{Theorem}[section]
\newtheorem{lemma}{Lemma}[section]

\newenvironment{proof}[1][Proof]{\textbf{#1.} }{\ \rule{0.5em}{0.5em}}
\usepackage{enumerate}
\immediate\write18{texcount -tex -sum  \jobname.tex > \jobname.wordcount.tex}
\usepackage[numbers]{natbib}
\usepackage{xcolor}
\usepackage{enumitem}
\usepackage{ulem}
\providecommand{\keywords}[1]
{
  \small	
  \textbf{\textit{Keywords---}} #1
}
\title{Bifurcation analysis of a prey-predator model with predator intra-specific interactions and ratio-dependent functional response}
\author{Claudio Arancibia--Ibarra$^{1,2}$, Pablo Aguirre$^{3}$, José D. Flores$^{4}$\\  \& Peter van Heijster$^{1,5}$\\
\small $^{1}$ School of Mathematical Sciences, Queensland University of Technology, Brisbane, Australia.\\
\small $^{2}$ Facultad de Ingeniería y Negocios, Universidad de Las Américas, Santiago, Chile.\\
\small $^{3}$ Departamento de Matemática, Universidad Técnica Federico Santa María, Valparaíso, Chile. \\
\small $^{4}$ Department of Computer Science, The University of South Dakota, South Dakota, USA. \\
\small $^{5}$ Biometris, Wageningen University and Research, Wageningen, Netherlands. \\
         \small Email: claudio.arancibia@hdr.qut.edu.au, pablo.aguirre@usm.cl, \\
         \small Jose.Flores@usd.edu, peter.vanheijster@wur.nl
}

\begin{document}
\maketitle

\begin{abstract}
We study the Bazykin predator-prey model with predator intraspecific interactions and ratio-dependent functional response and show the existence and stability of two interior equilibrium points. We prove that the model displays a wide range of different bifurcations, such as saddle-node bifurcations, Hopf bifurcations,  homoclinic bifurcations and Bogdanov-Takens bifurcations. We use numerical simulations to further illustrate the impact changing the predator per capita consumption rate has on the basin of attraction of the stable equilibrium points, as well as the impact of changing the efficiency with which predators convert consumed prey into new predators.
\end{abstract} \hspace{10pt}



\section{Introduction}
The original Lotka-Volterra predator-prey models~\cite{lotka} are relatively straightforward with simple functional forms for species’ growth and interactions. Empirical observations required successive changes to these assumptions, leading, inter alia, to the Bazykin models~\cite{turchin}. Temporal~\cite{aguirre3,greenhalgh,haque3,xiao} and spatio-temporal~\cite{banerjee3,bartumeus,shi,vazquez} ratio-dependent predator-prey models are becoming of more interest in ecology since these models better describe (when compared to the original models) both the theoretical and experimental predator-prey interactions~\cite{arditi,cosner}. Additionally, ratio-dependent predator-prey models are more suitable for predator-prey interactions when predators are seeking, capturing or killing other animals~\cite{arditi3,guti,kuang}. For instance, Jost and Arditi~\cite{jost2} found that ratio-dependent predator-prey models are more suitable for identifying the mathematical predator-prey functional response from real data. 

Gause type predator-prey models present an enrichment paradox~\cite{kuang,rosenzweig} which causes an increase in the equilibrium density of the predator but not in the prey, destabilising the community equilibrium. Gause type predator-prey models also present a biological control paradox~\cite{berryman2,luck} which refers to the fact that these type of models  cannot have a stable low prey density equilibrium. Kuang and Beretta~\cite{kuang} showed that ratio-dependent functional responses can resolve these paradoxes in Gause type predator-prey models.

Bazykin models~\cite{turchin,bazykin} with ratio-dependent functional response are examples of modified Gause type predator-prey models. The main difference between a standard Gause type model and a Bazykin model is that the latter considers the effect of adding a linear density dependence in the predator equation. The specific Bazykin model considered here is an extension of the model introduced by Arditi and Ginzburg~\cite{arditi} in which the authors showed that a ratio-dependent interaction is more realistic compared to a standard prey-dependent interaction. The Bazykin model is described by an autonomous two-dimensional system of ordinary differential equations, where the equations for the growth of the prey is a logistic-type function~\cite{turchin,aguirre1,flores2} and the growth of the predator is a function of the ratio of prey to predator abundance. The functional response is ratio-dependent, in which the consumption rate of prey depends on the numbers of predators and prey, i.e. $H(N,P)=qN/(N+aP)$~\cite{liu}. In other words, the ratio dependence is a type of predator dependence in which the functional response depends on the ratio of prey population size to predator population size, not on the absolute numbers of either species~\cite{abrams2}. In particular, it
is given by 
\begin{equation}\label{model1}
\begin{aligned}
\dfrac{dN}{dt} &=	 rN\left(1-\dfrac{N}{K} \right) - \dfrac{qNP}{N+aP}\,, \\
\dfrac{dP}{dt} &=	 \dfrac{cNP}{N+aP}-\mu_0P-\mu_1P^2 \,. 
\end{aligned}
\end{equation}
Here, $N(t)$ and $P(t)$ represent the proportion of the prey, respectively, predator population at time $t$; $r$ is the intrinsic growth rate of the prey; $K$ is the prey carrying capacity; $q$ is the per capita predation rate; $a$ is the amount of prey by which the predation effect is maximum; $c$ is the efficiency with which predators convert consumed prey into new predators; $\mu_0$ is the per capita death rate of predators and $\mu_1$ is predator death rate by density. 

Haque~\cite{haque2} studied a diffusive version of the Bazykin model \eqref{model1} and showed that, in absence of the diffusion, ratio-dependent predator-prey models are more appropriate for predator-prey interactions when the predators involve serious hunting processes. However, Haque focused on the case when system~\eqref{model1} has only one positive equilibrium point in the first quadrant (see Section~\ref{Num_eq}) and the author showed that there is always coexistence of both populations or the extinction of only the predator population. This manuscript extends the diffusive-free analysis studied by Haque~\cite{haque2} and its main aim is to study the bifurcation dynamics of~\eqref{model1}. In particular, we aim to understand the change in dynamics the intraspecific interactions and the ratio-dependent functional response causes. We will focus on the most general parameter setting for which system~\eqref{model1} can have up to two positive equilibrium points in the first quadrant and up to two limit cycles. As a result, system~\eqref{model1} supports complex dynamics such as the extinction, the coexistence, and the oscillation, of both populations over time. Furthermore, we will show that the model~\eqref{model1} supports different types of bifurcations such as saddle-node bifurcations, Hopf bifurcations and homoclinic bifurcations. 
We will also show the impact that
changing the predation rate has on the time series behaviour.

The basic properties of the Bazykin model \eqref{model1} are briefly described in Section~\ref{model}. In Section~\ref{result} we prove the stability of the equilibrium points and in Section~\ref{bif} we present the conditions for the different types of bifurcations. In addition, we discuss the impact changing the predation rate, i.e. $q$, or the efficiency with which predators convert consumed prey into new predators, i.e. $c$, has on the basins of attraction of the positive equilibrium points of system~\eqref{model1}. We further discuss the results and give the ecological implications in Section~\ref{con}.


\section{Preliminary Results}\label{model}
The ratio-dependent Bazykin predator-prey model is given by~\eqref{model1} and we only consider the model in the domain $\Omega=\{\left(N,P\right)\in\mathbb{R}^2,~N>0,~P>0\}$. In order to simplify the analysis we follow the nondimensionalisation approach of several other type of predator-prey models~\cite{arancibia7,arancibia4,arancibia5} and we introduce the dimensionless variables $(u,v,\tau)$ given by 
\begin{equation}\label{diff}
\begin{aligned}
& \varphi :\bar{\Omega}\times\mathbb{R}\rightarrow \Omega\times\mathbb{R}~\text{where}~\left(u,v,\tau\right)=\varphi\left(N,P,t\right) =\left(\dfrac{N}{K},\dfrac{aP}{K},\dfrac{rKt}{N+aP}\right).
\end{aligned}  
\end{equation}
Let $C:=c$, $M:=\mu_0/r$, $N:=\mu_1K/(ar)$ and $Q:=q/(ar)$, then~\eqref{model1} transforms into the nondimensionalised equivalent system
\begin{equation}\label{model2}
\begin{aligned}
\dfrac{du}{d\tau} & = u\left(1-u\right)\left(u+v\right)-Quv:=uW(u,v), \\ 
\dfrac{dv}{d\tau} & = Cuv-v\left(u+v\right)\left(M+Nv\right):=vR(u,v).
\end{aligned}  
\end{equation}
System~\eqref{model2} is defined in $\bar{\Omega}=\{(u,v)\in\mathbb{R}^2,~u\geq0,~v\geq0\}$ and $\varphi$~\eqref{diff} is a diffeomorphism for $(u,v)\neq(0,0)$ which preserve the orientation of time since $\det\left(\varphi\left(N,P,t\right)\right)=K^2\left(u^2+v^2\right)/(ar)>0$~\cite{andronov,chicone}. Moreover, system~\eqref{model2} is of Kolmogorov type~\cite{dumortier}, that is, the axes $u=0$ and $v=0$ are invariant. The $u$ nullclines are 
\begin{equation}\label{null_1}
u=0\quad\text{and}\quad v\left(Q-1+u\right)=u\left(1-u\right),
\end{equation}
 while the $v$ nullclines are 
 \begin{equation}\label{null_2}
 v=0\quad\text{and}\quad v=\dfrac{-\left(M+Nu\right)+\sqrt{\left(M+Nu\right)^2+4Nu\left(C-M\right)}}{2N}.
 \end{equation}

Hence, the equilibrium points on the axes for this system are $(0,0)$ and $(1,0)$\footnote{Note that the value $1$ of the equilibrium point $(1,0)$ relates to the rescaled prey carrying capacity $K$.}. In order to obtain the equilibrium point(s) inside the first quadrant, which we call positive interior equilibrium point(s), we rewrite equation~\eqref{null_1} as $v=u(1-u)/(Q-1+u)$ (assuming $Q-1+u\neq0$).
Then, replacing $v$ into $R(u,v)=0$~\eqref{model2} we obtain: 
\begin{equation}\label{Null_int}
\begin{aligned}
\left(C+NQ\right)u^2+\Sigma_1 u+\left(Q-1\right)\Sigma_2=0
\end{aligned}
\end{equation}
with 
\begin{equation}\label{S12}
\begin{aligned}
\Sigma_1&=2C\left(Q-1\right)-Q\left(M+N\right) = 2 \Sigma_2 + Q(M-N),~\text{and}\\
\Sigma_2&=C\left(Q-1\right)-MQ = Q(C-M)-C. 
\end{aligned}
\end{equation}
Similarly, we can rewrite equation~\eqref{null_2} as $u=v(M+Nv)/(C-M-Nv)$ (assuming $C-M-Nv\neq0$). Then, replacing $u$ into $W(u,v)=0$~\eqref{model2} we obtain: 
\begin{equation}\label{Null_int_v}
N\left(C+NQ\right)v^2+\Sigma_3 v+\left(C-M\right)\Sigma_2=0
\end{equation}
with 
\begin{equation}\label{S34}
\Sigma_3=2NQ(M-C)+C(M+N) = -2N \Sigma_2 + C(M-N)= -N\Sigma_1+(QN+C)(M-N)\,.
\end{equation}

Therefore, the interior equilibrium point(s) $P_{1,2}=(u_{1,2},v_{1,2})$ of system~\eqref{model2} are determined by the solution(s) of the equations~\eqref{Null_int} and~\eqref{Null_int_v} which are given by: 
\begin{equation}\label{epp1p2}
u_{1,2} = \dfrac{-\Sigma_1 \pm Q\sqrt{\Delta}}{2\left(C+NQ\right)}~~~\text{and}~~~v_{1,2} = \dfrac{-\Sigma_3 \pm C\sqrt{\Delta}}{2N\left(C+NQ\right)}.\\
\end{equation} 

Here, 
\begin{equation}\label{delta}
\Delta=\left(M-N\right)^2-4N\left(C\left(Q-1\right)-MQ\right)=\left(M-N\right)^2-4N\Sigma_2,
\end{equation}
and $u_{1}<u_{2}<1$ and $v_{1}<v_{2}$. 

\subsection{Number of positive equilibrium points}\label{Num_eq}

We describe the different configurations for the solutions of equations~\eqref{Null_int} and~\eqref{Null_int_v} given in \eqref{epp1p2}. In particular, we determine the number of positive interior equilibrium points. 
	
\begin{figure}
\centering
\includegraphics[width=10cm]{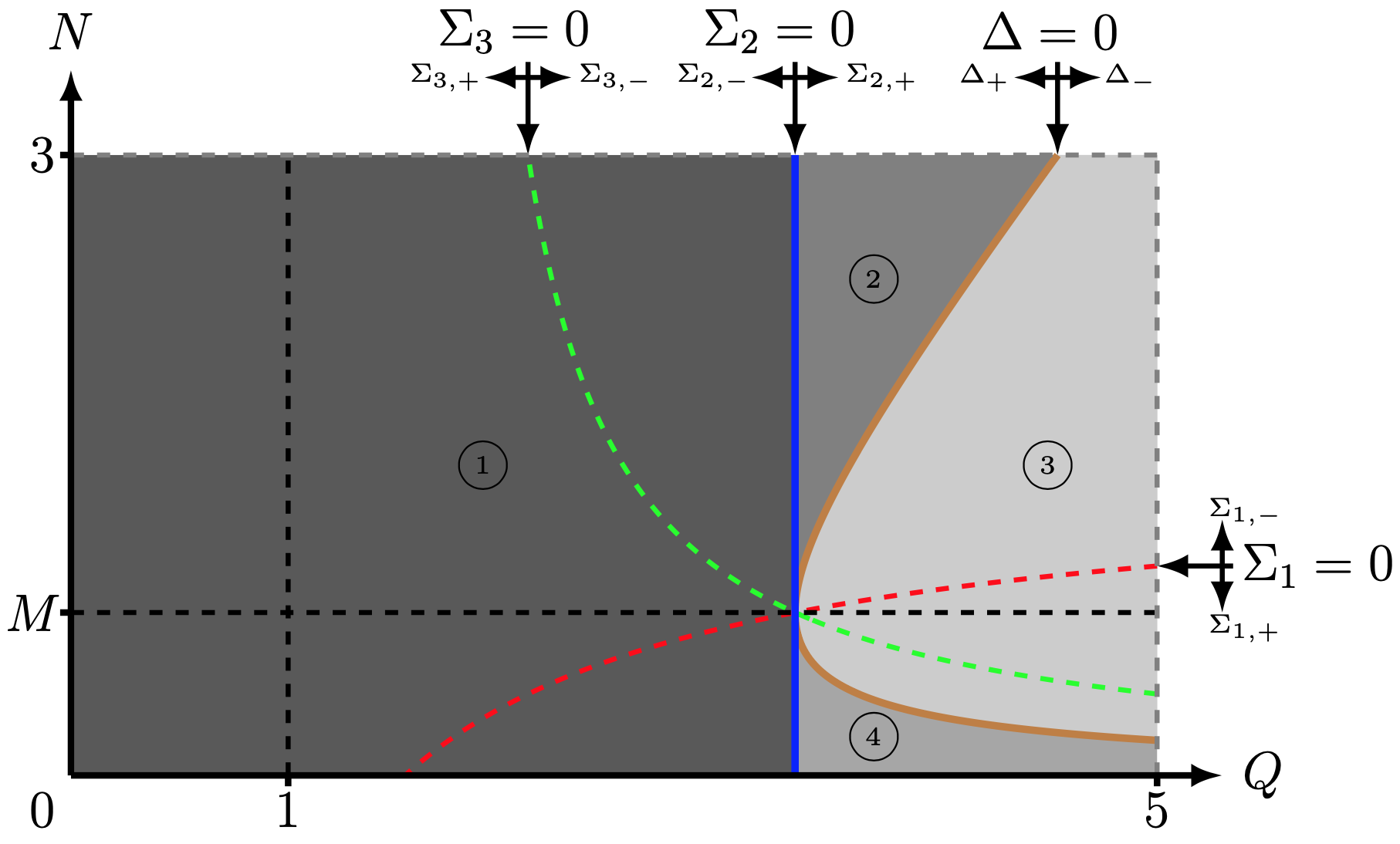}
\caption{Parametric diagram with the different conditions for the configuration of $u_{1,2}$ and $v_{1,2}$~\eqref{epp1p2} in the $(Q,N)$-plane for $C=1.5$ and $M=1.05$ fixed, such that $C>M$.  In region 1, $\Sigma_2:=C\left(Q-1\right)-MQ<0$ and system~\eqref{model2} has one positive interior equilibrium point (case~\ref{i}).
In region 3, $\Delta:=\left(M-N\right)^2-4N\Sigma_2<0$  and system~\eqref{model2} has no positive interior equilibrium points (case~\ref{DEL}).
In region 2, where $\Sigma_2>0,\Delta>0$ and $N>M$, system~\eqref{model2} has two positive equilibrium points (case~\ref{AA}), while 
there are no positive interior equilibrium points in region 4 where $\Sigma_2>0,\Delta>0$ and $N<M$ (case~\ref{BB}).}
\label{Fig1}
\end{figure}

\begin{enumerate}
\item \label{D} If $C\leq M$, then $u_{i}$ and $v_{i}$, $i=1,2$, cannot be both positive since $u=v(M+Nv)/(C-M-Nv)$. Therefore, there are no positive interior equilibrium points.
\item \label{C} If $C>M$ and
\begin{enumerate}
 \item \label{i} $\Sigma_2<0$, then $\Delta>0$ and $(C-M)\Sigma_2<0$. 
Hence, the two solutions of Equation \eqref{Null_int_v} are real, but only $v_2$ is positive. 
Equation~\eqref{Null_int} also has two real solutions and $u_2$ is always positive since 
\begin{enumerate}
\item for $Q>1$ we have that $(Q-1)\Sigma_2<0$; and
\item for $Q \leq 1$ we have that $\Sigma_1<0$. 
\end{enumerate}
Therefore, system~\eqref{model2} has one positive interior equilibrium point $P_{2} = (u_{2}, v_{2})$ in this case. See region 1 in Figure~\ref{Fig1}. 
\item \label{ii} $\Sigma_2=0$, then one of the equilibrium points $(u_{1,2},v_{1,2})$ merges with $(0,0)$. The remaining equilibrium point is in the interior of the first quadrant if both $\Sigma_1$ and $\Sigma_3$ are negative. Since $\Sigma_2=0$, both of these inequalities simplify to $N>M$. Therefore, in this case system~\eqref{model2} has one positive interior equilibrium point $P_{2} = (u_{2}, v_{2})$ if $N>M$ and no positive interior equilibrium points if $N \leq M$.  

\item \label{SG} $\Sigma_{2}>0$ and 
\begin{enumerate}
\item \label{DEL} $\Delta<0$, then then system~\eqref{model2} has no positive equilibrium points in the first quadrant, see region 3 in Figure~\ref{Fig1}.
\item \label{DEL0} $\Delta=0$, then both roots merge. So, $u_1=u_2$ and $v_1=v_2$. This root of order two is in the first quadrant if and only if both $\Sigma_1$ and $\Sigma_3$ are negative.
Along $\Delta=0$ the expression for $2N\Sigma_1$ becomes $$ 2N \Sigma_1 = (M-N)\left( M+N(2Q-1)\right).$$ 
So, since $\Sigma_2>0$ implies that $Q>1$, we have that $\Sigma_1>0$ along $\Delta=0$ if $M>N$ and $\Sigma_1<0$ along $\Delta=0$ if $M<N$.
Similarly, the expression for $2\Sigma_3$ along $\Delta=0$  is $$ 2 \Sigma_3 = (M-N)\left( 2C+N-M\right).$$ 
So, since $C>M$, we have that $\Sigma_3>0$ along $\Delta=0$ if $M>N$ and $\Sigma_3<0$ along $\Delta=0$ if $M<N$.
Therefore, in this case system~\eqref{model2} has one positive interior equilibrium point $P_{2} = (u_{2}, v_{2})$ of order two if $N>M$ and no positive interior equilibrium points if $N \leq M$.   
\item \label{DEL2} $\Delta>0$ and
\begin{enumerate}
\item \label{AA}
$N>M$, then $\Sigma_1<0$ -- which follows directly from the above computation in case~\ref{DEL0} upon replacing the =-sign by the <-sign -- and $\Sigma_3=-2N \Sigma_2 + C(M-N)<0$.
Thus, 
since $\Sigma_2>0$ implies that $Q>1$, we have that both the $u$-intercept and $v$-intercept are positive. Consequently, system~\eqref{model2} has two positive interior equilibrium points. See region 2 in Figure~\ref{Fig1}.
\item \label{BB} 
$N \leq M$, then $\Sigma_1= 2\Sigma_2+Q(M-N) \geq 0$ and
$\Sigma_3 \geq 0$ -- which follows directly from the above computation in case~\ref{DEL0} upon replacing the =-sign by the $\geq$-sign. Furthermore, we still have that both the $u$-intercept and $v$-intercept are positive. Consequently, system~\eqref{model2} now has no positive interior equilibrium points. See region 4 in Figure~\ref{Fig1}.
\end{enumerate}
\end{enumerate}
\end{enumerate}
\end{enumerate}

\section{Main Results}\label{result}
In this section, we discuss the stability of the equilibrium points of system~\eqref{model2} for $Q>1$ which represents the most general case of the number of equilibrium points. 
Since $Q:=q/(ar)$, the condition $Q>1$ corresponds to the assumption that the predation rate is larger than the product 
between the prey's intrinsic growth rate and the number of prey necessary for getting the maximum predation effect. 
Note that the case when $Q=1$ and $C>M$ (a subcase of case \ref{i} in Section~\ref{Num_eq}) was studied by Haque~\cite{haque2}.
\begin{lemma}\label{bou}
The set $\Gamma=\{\left(u,v\right)\in\bar{\Omega},~0\leq u\leq1,~v\geq0\}$ is an invariant region and all solutions of~\eqref{model2} which are initiated in the first quadrant are bounded and eventually end up in $\Gamma$.
\end{lemma}
\begin{proof}
This proof follows the proofs of Arancibia \textit{et al}.~\cite{arancibia4,arancibia3}. Since the system~\eqref{model2} is of Kolmogorov type the coordinates axes are invariant~\cite{freedman}. Moreover, if $v=0$ then $du/d\tau=u^2(1-u)$ and if $u=0$ then $dv/d\tau=-v^2(M+Nv)$. So, any trajectory with initial point on the positive vertical $v$-axis tends to zero and any trajectory with initial point on the positive horizontal $u$-axis tend to $u=1$. Next, setting $u=1$ in the first equation of system~\eqref{model2}, we have $du/dt=-Qv<0$ and setting $u>1$, we also have $du/dt<0$ since $1-u<0$. Thus for any initial condition initiated in the first quadrant, the corresponding trajectory enters and remains in $\Gamma$ independent of the sign of $dv/dt$.

To finalise the proof we show that no trajectory in the open region $\Gamma$ converges to infinity as $\tau \to \infty$. To show that solutions are bounded it is enough to find a $v^*$ such that $dv/d\tau<0$ for $v \geq v^*$ and $0\leq u\leq1$. The equation $dv/d\tau$ can be written as $dv/d\tau=-Cv^3+\mathcal{O}(v^2,uv^2)$ which implies that we can always find such as $v^*$. Therefore, all trajectories end up in $\Gamma_1=\{\left(u,v\right)\in\bar{\Omega},~0\leq u\leq1,~0\leq v\leq v^*\}$. 
\hfill\end{proof}

\subsection{The nature of the equilibrium points}
To determine the nature of the equilibrium points we compute the Jacobian matrix $J(u,v)$ of~\eqref{model2}
\begin{equation}\label{jac}
J(u,v)=\begin{pmatrix}
 2u+v-Qv-2uv-3u^2 & -u(Q+u-1) \\
-v(M-C+Nv) & Cu-3Nv^2-Mu-2Mv-2Nuv
\end{pmatrix}\,.
\end{equation}
\subsubsection{Equilibrium points on the axes}
\begin{lemma}\label{1,0}
The equilibrium point $(1,0)$ is a saddle if $C>M$ and a stable node if $C<M$.
\end{lemma}
\begin{proof}
The result follows direct from the Jacobian matrix~\eqref{jac} evaluated at $(1,0)$
\begin{equation*}
J(1,0)=\begin{pmatrix}
 -1 & -Q \\
0 & C-M
\end{pmatrix}
\end{equation*} 
combined with the Hartman-Grobman theorem.
\hfill\end{proof}

\begin{figure}
\centering
\includegraphics[width=12.8cm]{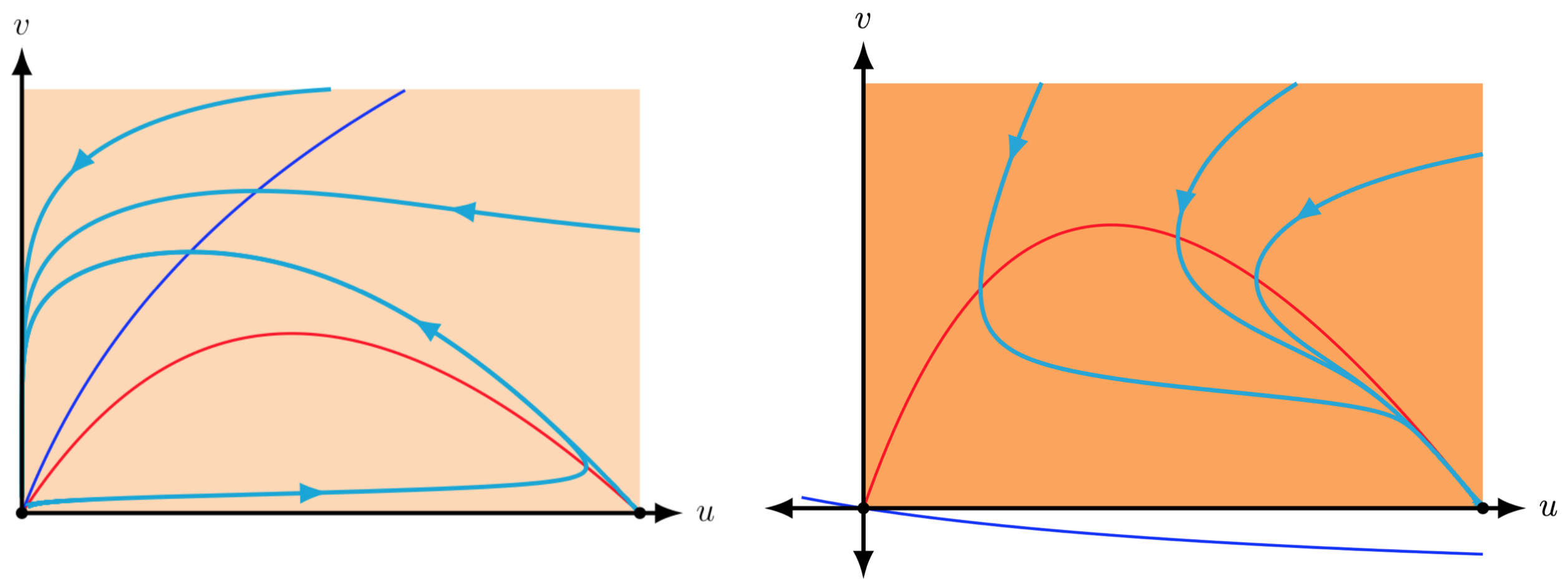}
\caption{The blue (red) curve represents the prey (predator) nullcline, the light orange region represents the basin of attraction of $(0,0)$ and the dark orange region represents the basin of attraction of $(1,0)$. In the left panel $Q=3.05$, $C=10.05$, $M=1.05$, and $N=10.0$ such that we are in case~\ref{DEL} of Section~\ref{Num_eq}. The two nullclines do not intersect and there are no positive equilibrium points. The equilibrium point $(1,0)$ is a saddle since $C>M$, see Lemma~\ref{1,0}, and the origin $(0,0)$ is an attractor with an elliptical sector since both $C>M+1$ and $Q>M+1$, see Theorem~\ref{00}.  In the right panel $Q=1.8$, $C=0.205$, $M=0.22$, and $N=0.25$ such that we are in case~\ref{D} of Section~\ref{Num_eq}. The two nullclines do not intersect and there are no positive equilibrium points. The origin $(0,0)$ is a saddle point and, since $C<M$, $(1,0)$ is an attractor for positive initial conditions.}
\label{Fig2}
\end{figure}

In other words, if $C<M$  then $(1,0)$ is an attractor for positive initial conditions, see the right panel of Figure~\ref{Fig2}.

\begin{figure}
\centering
\includegraphics[width=9cm]{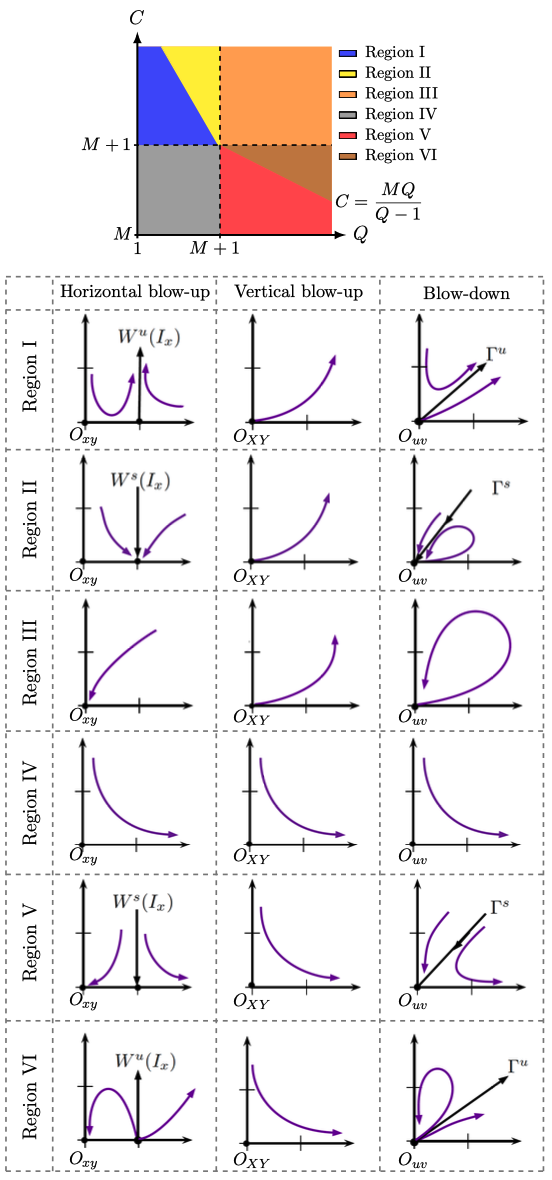}
\caption{Parametric diagram with the different structures in a neighbourhood of the origin $(0,0)$ in the $(Q,C)$-parameter-space for system~\eqref{model2}. The first column shows the horizontal blow up, the second column the vertical blow up, while the third column shows the blow down.}
\label{Fig3}
\end{figure}

To understand the dynamics near the origin we assume, 
in addition to the assumption $Q>1$, that $C>M$. That is, we assume that the equilibrium point $(1,0)$ is a saddle point, see Lemma~\ref{1,0}. In other words,  we assume that the efficiency with which predators convert consumed prey into new predators is bigger that the ratio between the per capita death rate of predators and the prey intrinsic growth rate.
We divide our $(Q,C)$-parameter space in six regions
\begin{itemize}
\item Region I: $M+1<C<MQ/(Q-1)$;
\item Region II: $C>MQ/(Q-1)$ and $1<Q<M+1$;
\item Region III: $C>M+1$ and $Q>M+1$;
\item Region IV: $M<C<M+1$ and $1<Q<M+1$;
\item Region V: $M<C<MQ/(Q-1)$ and $Q>M+1$; and
\item Region VI: $MQ/(Q-1)<C<M+1$,
\end{itemize} 
see also Figure~\ref{Fig3}. 
\begin{theorem}\label{00}
If we assume that $Q>1$ and $C>M$, then the origin $(0,0)$ in system~\eqref{model2} is a non-hyperbolic degenerate point~\cite{berezovskaya1,berezovskaya2}. Moreover, a neighbourhood of the origin $(0,0)$ presents six types of topologically different structures in the first quadrant of the phase plane:
\begin{itemize}
\item a saddle sector and a repelling sector in Region I; 
\item an attracting sector and elliptic sector in Region II;
\item an elliptic sector in Region III;
\item a saddle sector in Region IV;  
\item an attracting sector and a saddle sector in Regions V; and
\item an elliptic sector and a repelling sector in Region VI.
\end{itemize}    
\end{theorem}
We refer to Figure~\ref{Fig3} for exemplary sketches of the dynamics near the origin in the six different regions.
\begin{proof}
First we observe that setting $u=0$ in system~\eqref{model2} gives $dv/d\tau=-v^2(M+Nv)\leq0$ for $v\geq0$. That is, any trajectory starting along the $v$-axis converges
to the origin $(0,0)$. Also, setting $v=0$ in system~\eqref{model2} gives $du/d\tau=u^2(1-u)$ and any orbits starting along the $u$-axis near the origin converges to the equilibrium point $(1,0)$. Next, to analyse the dynamics in a neighbourhood of the origin, we consider the vertical blow-up~\cite{dumortier} given by the transformation 
\begin{equation}
(u,v)\to \left(xy,y\right)~\text{and the time rescaling}~\tau\to \dfrac{t}{y}. \label{eq1}
\end{equation}
The transformation~\eqref{eq1} \textit{blows up} the origin of system~\eqref{model2} in the entire $x$-axis~\cite{dumortier}. The goal is to analyse the equilibrium points on the nonnegative half axis $x\geq 0$, $y=0$, in the new coordinate system
\begin{equation}\label{eq2}
\begin{aligned}
\dfrac{dx}{dt} & = x((x+y)y(N-x)+(x+1)(1+M)-Q-Cx), \\
\dfrac{dy}{dt} & = y(Cx-(x+1)(M+Ny)).
\end{aligned}
\end{equation}
System~\eqref{eq2} has up to two equilibrium points on the nonnegative horizontal $x$-axis: the origin $O_{xy}$ and, for $(C,Q)$ in regions I, II, V or VI, a second equilibrium point $I_x=(\mu,0)$ with $\mu=(1+M-Q)/(C-M-1)$. 
The corresponding Jacobian matrix at $O_{xy}$ is
\begin{equation*}
J(O_{xy})=\begin{pmatrix}
 1+M-Q && 0 \\
0 && -M
\end{pmatrix}
\end{equation*}
with eigenvalues 
\begin{equation*}\lambda_1(O_{xy})=1+M-Q~\text{and}~\lambda_2(O_{xy})=-M.\end{equation*}
Hence, $O_{xy}$ is
\begin{itemize}
\item a saddle if $1<Q<M+1$, see Regions I, II and IV in Figure~\ref{Fig3}; and
\item a stable node if $Q>M+1$, see Regions III, V and VI in Figure~\ref{Fig3}.
\end{itemize}

At $I_x$ we get the Jacobian matrix
\begin{equation*}
J(I_x)=\begin{pmatrix}
 -(1+M-Q) && \dfrac{(Q-C)(N(C-1)-M(1+N)+Q-1)(1+M-Q)}{(1+M-C)^3} \\
0 && \dfrac{\Sigma_2}{(1+M-C)}
\end{pmatrix}
\end{equation*}
with eigenvalues 
\begin{equation}\label{ev_Ix}
\lambda_1(I_{x})=-1-M+Q~\text{and}~\lambda_2(I_{x})=\dfrac{\Sigma_2}{1+M-C}.
\end{equation}
Hence, $I_x$ is
\begin{itemize}
\item  a saddle in regions I and V, see Figure~\ref{Fig3}.
\item an stable node in region II, see Figure~\ref{Fig3}.
\item an unstable node in region VI, see Figure~\ref{Fig3}.
\end{itemize}

We repeat the blow-up procedure to analyse the behaviour of system~\eqref{model2} near the $v$-axis. Thus we consider the horizontal blow-up~\cite{dumortier} given by the transformation
\begin{equation}
(u,v)\to \left(X,XY\right)~\text{and the time rescaling}~\tau\to\dfrac{t}{X} \label{eq3}
\end{equation}
The goal now is to analyse the equilibrium point at the nonnegative half axis $X=0$, $Y\geq0$, in the new system coordinates which is given by setting~\eqref{eq3} in system~\eqref{model2}
\begin{equation}\label{eq4}
\begin{aligned}
\dfrac{dX}{dt } & = X(1-X+Y(1-Q-X)), \\
\dfrac{dY}{dt } & = Y(C-1+Y(Q-1)-(Y+1)(M-X+NXY)).
\end{aligned}
\end{equation}
System~\eqref{eq4} has again up to two equilibrium points in the nonnegative vertical $Y$-axis. The equilibrium point $(0,0)$ and, for $(C,Q)$ in regions I, II, V or VI, a second equilibrium point $I_Y=(0,\mu^*)$ with $\mu^*=(C-M-1)/(1+M-Q)$. 

At the origin $O_{XY}=(0,0)$ we get the Jacobian matrix 
\begin{equation*}
J_{(O_{XY})}=\begin{pmatrix}
 1 & 0 \\
0 & C-1-M
\end{pmatrix}
\end{equation*}
with eigenvalues \begin{equation*}\lambda_1(O_{XY})=1~\text{and}~\lambda_2(O_{XY})=C-1-M.\end{equation*}
Hence, the stability of $O_{XY}$ depends on the parameters $C$ and $M$. In particular, $O_{XY}$ is
\begin{itemize}
\item a saddle if $M<C<M+1$, see Regions IV, V and VI in Figure~\ref{Fig3}.
\item an unstable node if $C>M+1$, see Regions I, II and III in Figure~\ref{Fig3}.
\end{itemize}

The Jacobian matrix at $I_Y$ is given by
\begin{equation*}
J(I_Y)=
\begin{pmatrix}
\dfrac{-\Sigma_2}{M-Q+1} && 0\\
\dfrac{(C-Q)(C-M-1)(M+N-Q-CN+MN+1)}{(M-Q+1)^3} && 1+M-C
\end{pmatrix}
\end{equation*}
with eigenvalues 
\begin{equation}\label{ev_Iy}
\lambda_1(I_{Y})=\dfrac{-\Sigma_2}{M-Q+1}~\text{and}~\lambda_2(I_{Y})=1+M-C.
\end{equation}
Upon comparing~\eqref{ev_Ix} with~\eqref{ev_Iy} we note that the stability of $I_Y$ is the same as the stability of $I_x$. 

Taking the inverse of~\eqref{eq1} and~\eqref{eq3} the equilibrium points $I_x$ and $I_Y$ collapse to the origin $(0,0)$ of~\eqref{model2}, and the stable and unstable eigenvectors are mapped to the curves $\Gamma^s$ and $\Gamma^u$ respectively, see the blow-down pictures in the right column of Figure~\ref{Fig3}. Therefore, it follows that the origin is a non-hyperbolic degenerate point with the properties as stated in Theorem~\ref{00}. 
\hfill\end{proof}

\subsubsection{The positive equilibrium points}
Next, we consider the stability of the two positive equilibrium points $P_{1,2}$ of system~\eqref{model2} in the interior of $\Gamma$. These equilibrium points lie on the curve $v=u\left(1-u\right)/\left(Q-1+u\right)$ 
such that $W(u,v)=0$ \eqref{model2}, and they only exist if the system parameters are such that $\Delta>0$~\eqref{delta}. 
The Jacobian matrix of system~\eqref{model2} at these equilibrium points
becomes
\begin{equation}\label{jac_2}
J(u,v)=\begin{pmatrix}
 u(1-2u-v) & u(1-Q-u) \\
v(C-M-Nv) & -v(M+N(u+2v))
\end{pmatrix} \,.
\end{equation}
The determinant and the trace of the Jacobian matrix~\eqref{jac_2} are given by
\begin{equation}\label{det}
\begin{aligned}
&\det\left(J(u,v)\right)=CN(Q-1)(u+v)(3(u+v)-2)+M(2(u+v)-Q)~\text{and}\\
&\tr\left(J(u,v)\right)=u(1-2u-v)-v(M+N(u+2v))=:T(u,v). 
\end{aligned}
\end{equation}
This gives the following results.

\begin{theorem}\label{p1}
Let the system parameters of~\eqref{model2} be such that $C>M$, $\Sigma_{2}>0$, $\Delta>0$~\eqref{delta} and $N>M$ (case~\ref{AA} of Section~\ref{Num_eq}), then there are two interior equilibrium points and the equilibrium point $P_1=\left(u_1,v_1\right)$ defined in~\eqref{epp1p2} is a saddle point.
\end{theorem}
\begin{proof}
Since $P_1$ is in the first quadrant, we know that $0<u_1+v_1=\left(N-M-\sqrt{\Delta}\right)/\left(2N\right)$, with $\Delta$ defined in~\eqref{delta}. Evaluating the determinant~\eqref{det} at $P_1=\left(u_1,v_1\right)$ gives 
\begin{equation*}
\det\left(J\left(u_1,v_1\right)\right)=-\dfrac{\sqrt{\Delta}\left(N-M-\sqrt{\Delta}\right)}{2N}<0.
\end{equation*} 
Hence, the equilibrium point $P_1$ is a saddle point. \hfill\end{proof}

\begin{theorem}\label{p2}
Let the system parameters of~\eqref{model2} be such that we are in cases~\ref{AA} or~\ref{i} of Section~\ref{Num_eq}, then system ~\eqref{model2} has at least one positive interior equilibrium point $P_2=\left(u_2,v_2\right)$. Therefore, the stability of equilibrium point $P_2$ defined in~\eqref{epp1p2} is:
\begin{enumerate}[label=(\roman*)]
\item stable if $T\left(u_2,v_2\right)<0$;
\item unstable if $T\left(u_2,v_2\right)>0$; and
\item a weak-focus if $T\left(u_2,v_2\right)= 0$,
\end{enumerate}
where $T$ is defined in \eqref{det}.
\end{theorem}
\begin{proof}
Evaluating the determinant~\eqref{det} at $P_2=\left(u_2,v_2\right)$ gives
\begin{equation*}
\det\left(J\left(u_2,v_2\right)\right)=\dfrac{\sqrt{\Delta}\left(N-M+\sqrt{\Delta}\right)}{2N}>0,
\end{equation*} 
since $u_2+v_2=(N-M+\sqrt{\Delta})/(2N)>0.$
The results now follow by analysing the trace~\eqref{det} of the Jacobian matrix~\eqref{jac_2} evaluated at $P_2$. The function $T\left(u,v\right)$~\eqref{det} at the equilibrium point $P_2$ is given by
\begin{equation*}
T\left(u_2,v_2\right)=\dfrac{1}{N\left(C+NQ\right)^2}\left(T_1\sqrt{\Delta}+T_2\right)\,,
\end{equation*}
with
\begin{equation*}
\begin{aligned}
T_1\big(C,N,M,Q\big)=&NQ^2\big(3C-M-N-CN^2+MN^2\big)+CQ\big(C-3N-CN^2-2CN^3+MN^2+2MN^3\big) \\
&+C^2\big(MN^2+N^3-1\big)
\end{aligned}
\end{equation*}
and
\begin{equation*}
\begin{aligned}
T_2\big(C,N,M,Q\big)=&4N^2 Q^3\big(C-M\big)-N Q^2\big(CN^3-MN^3-2C^2N^2+2C^2N^3-M^2N^2+2M^2N^3+CM\\
&+3CN-2MN+M^2+N^2+3CMN^2-4CMN^3\big)+CQ\big(-3CN^3+2CN^4+3MN^3\\
&-2MN^4+2C^2N^2+2C^2N^3+M^2N^2-2M^2N^3-CM+CN+MN-N^2-3CMN^2\big)\\
&-C^2\big(-M+N+2CN^2+2CN^3-2MN^2+M^2N^2+N^4\big)\,.
\end{aligned}
\end{equation*}
Therefore, the sign of the trace $T\left(u_2,v_2\right)$~\eqref{det}, and thus the behaviour of $P_2$, depends on the parity of $T_1\sqrt{\Delta}+T_2$. Evaluating $T_1$, $T_2$ and $\Delta$ at $(C,M,N,Q)=(0.363,0.16,0.25,1.6)$ gives $T_1\sqrt{\Delta}+T_2=-0.036091<0$, while $T_1\sqrt{\Delta}+T_2=0.025983>0$ at $(C,M,N,Q)=(0.363,0.16,0.25,1.8)$. By the continuous dependence of the trace on the parameter $Q$, we get that there exist $Q\in(1.6,1.8)$ for which the trace must be zero. 
\hfill\end{proof}

\begin{figure}
\centering
\includegraphics[width=5.3cm]{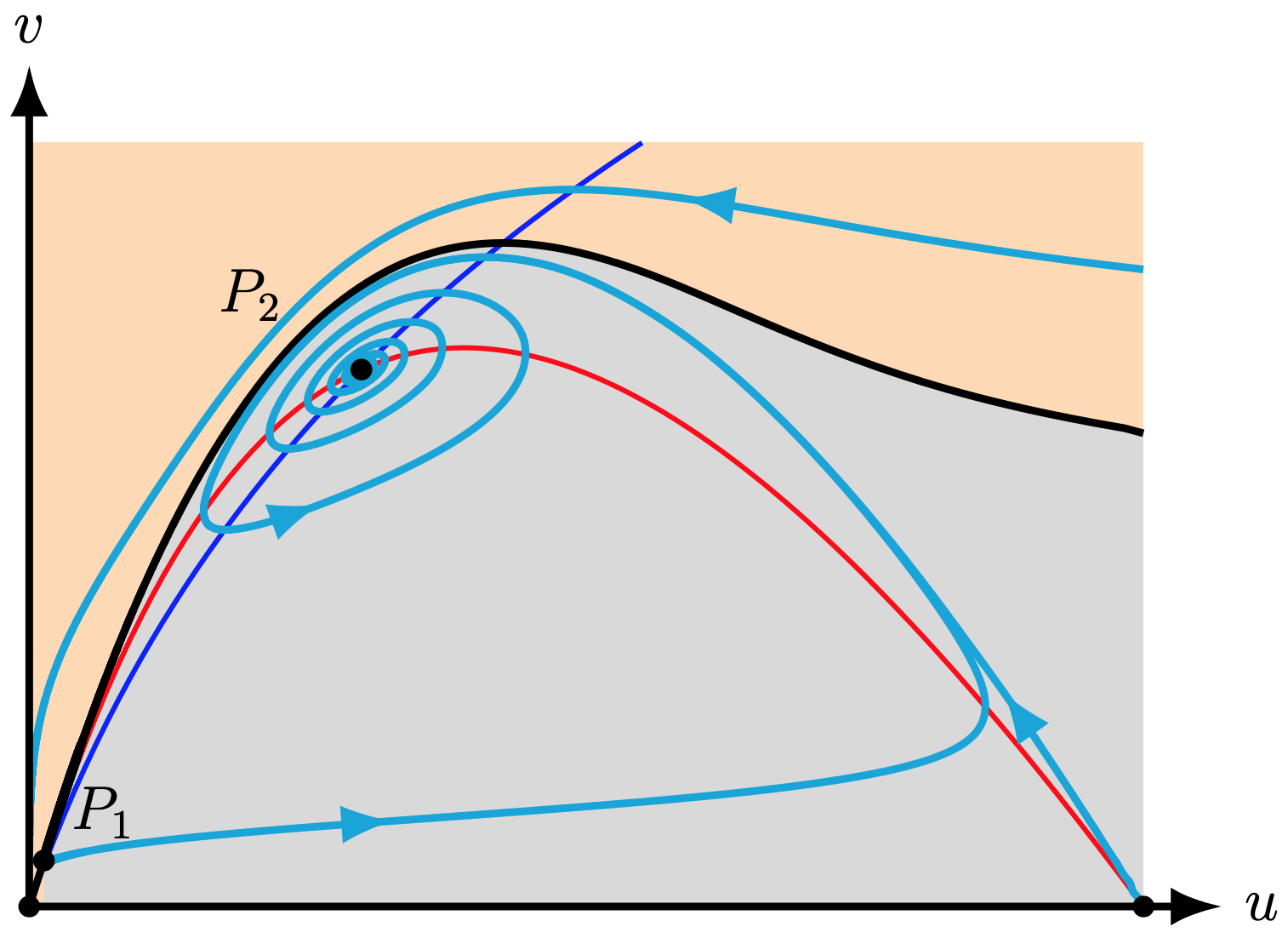}
\includegraphics[width=5.3cm]{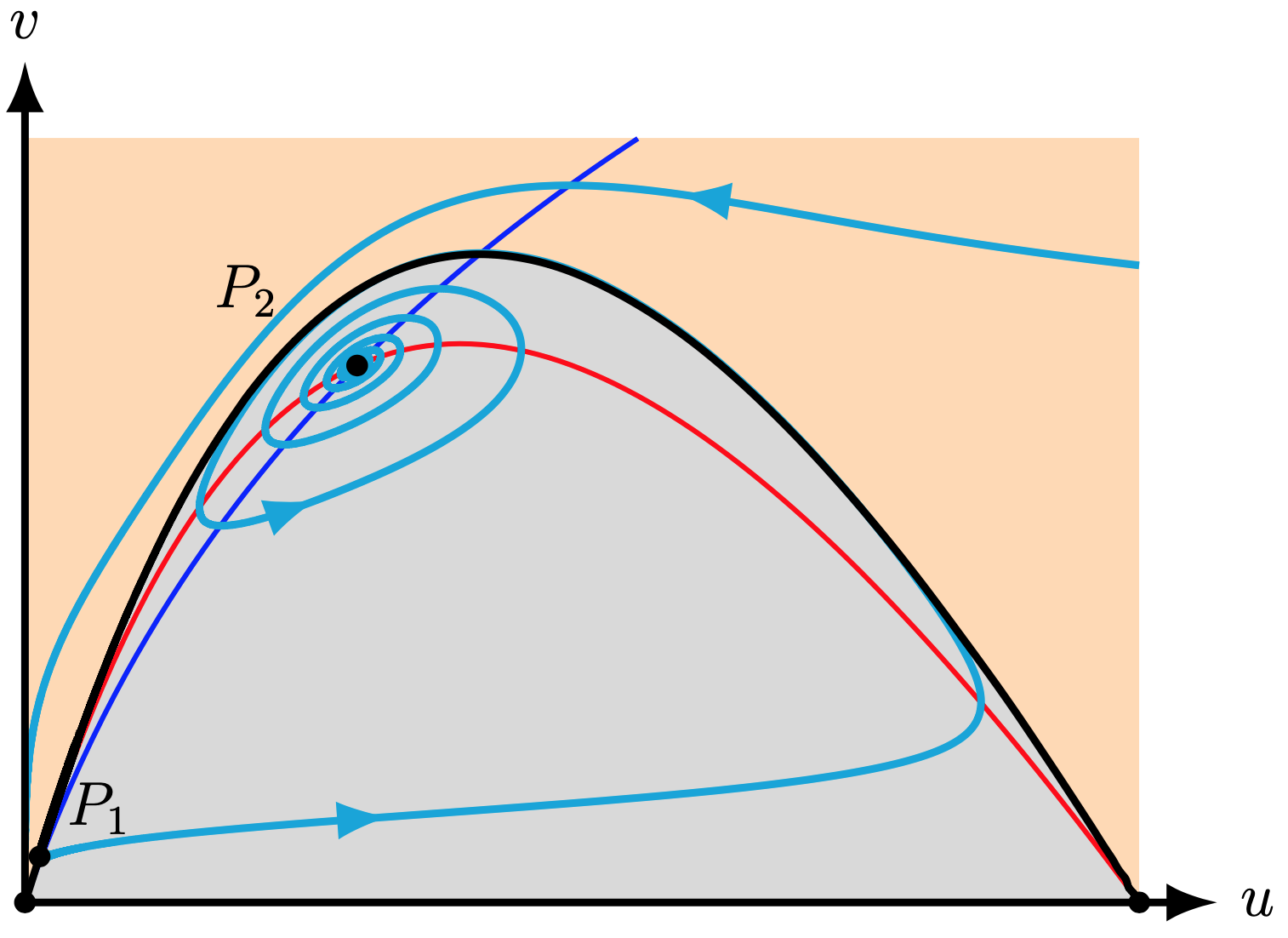}
\includegraphics[width=5.3cm]{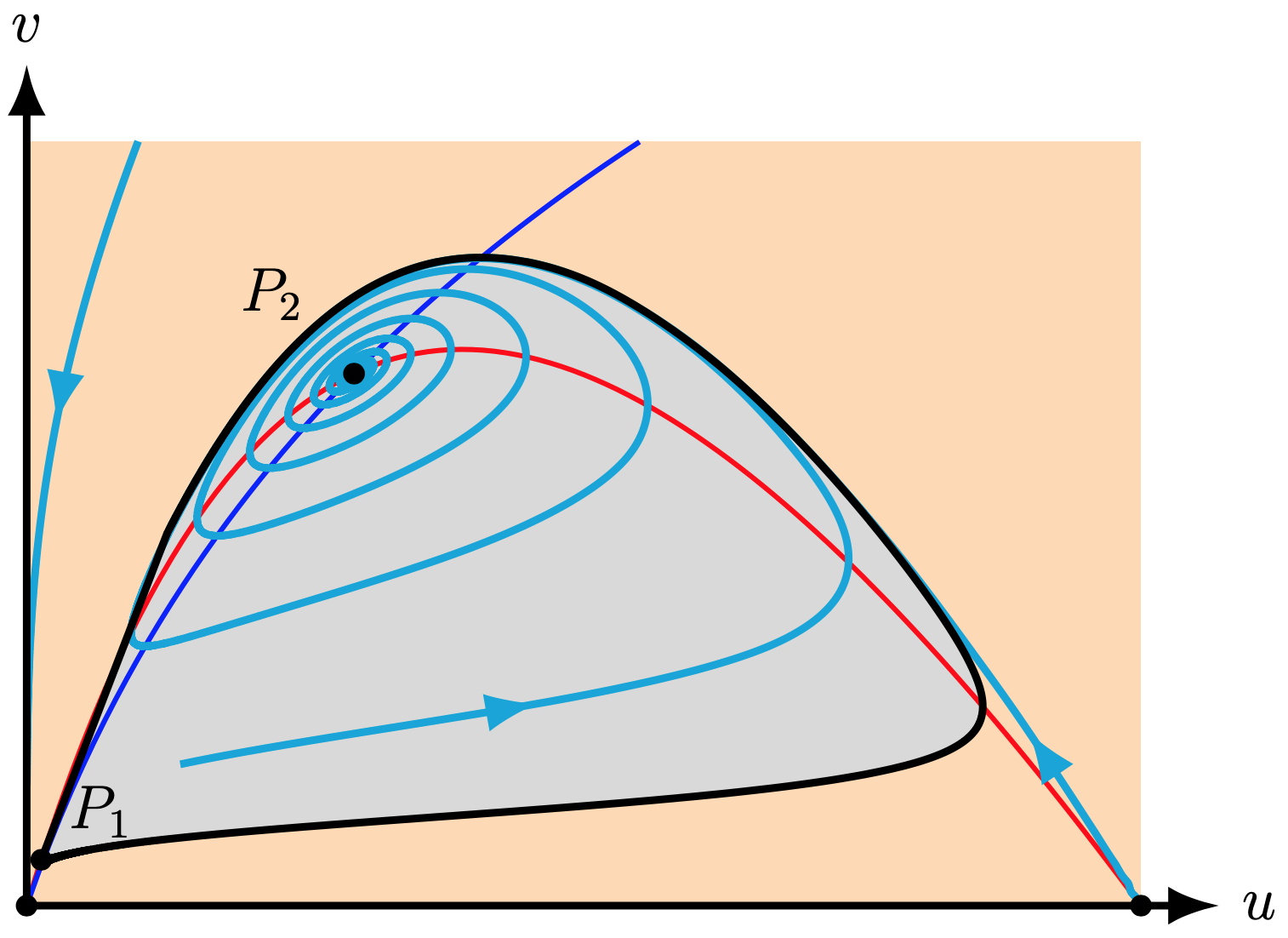}
\includegraphics[width=5.3cm]{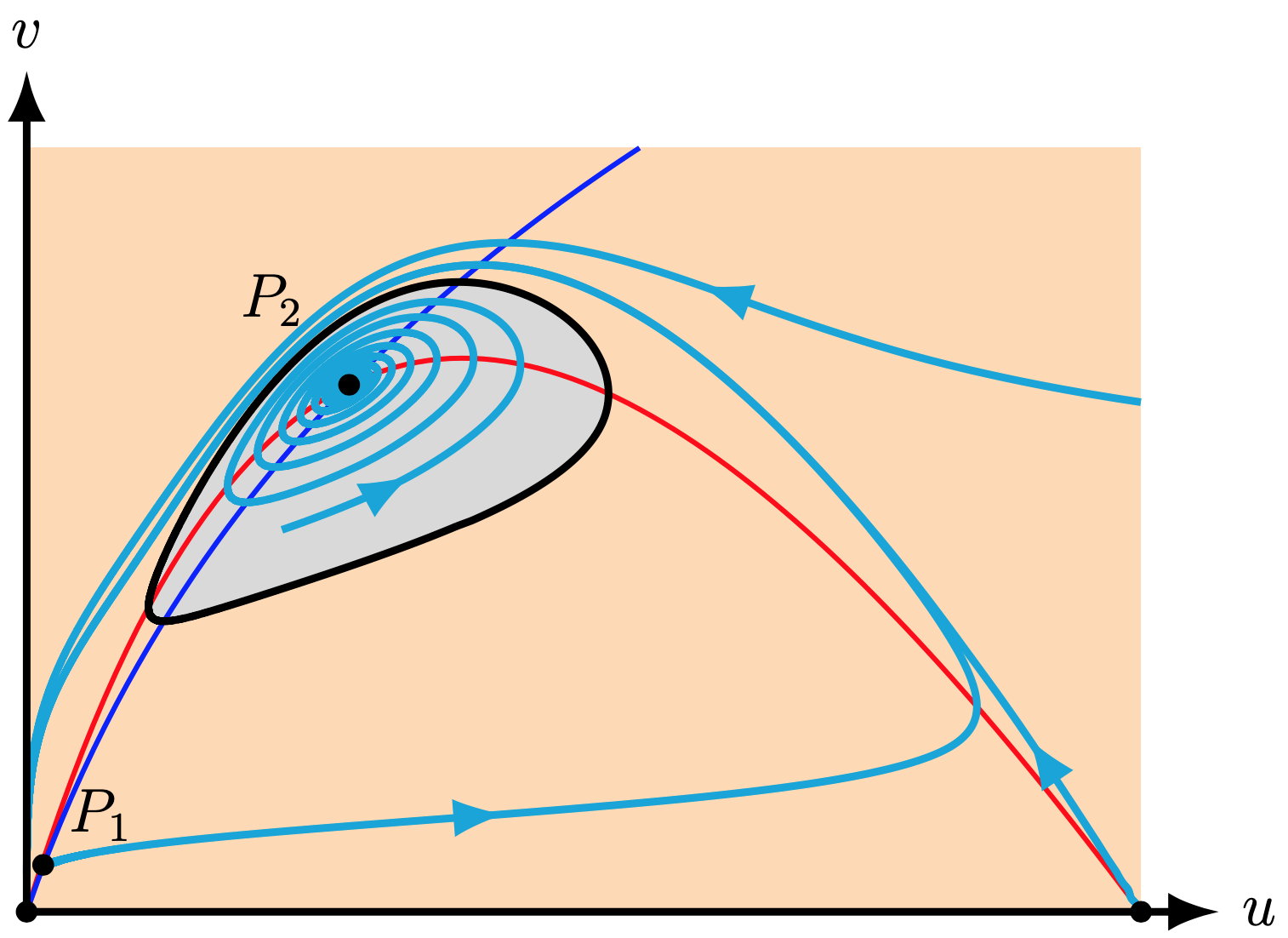}
\includegraphics[width=5.3cm]{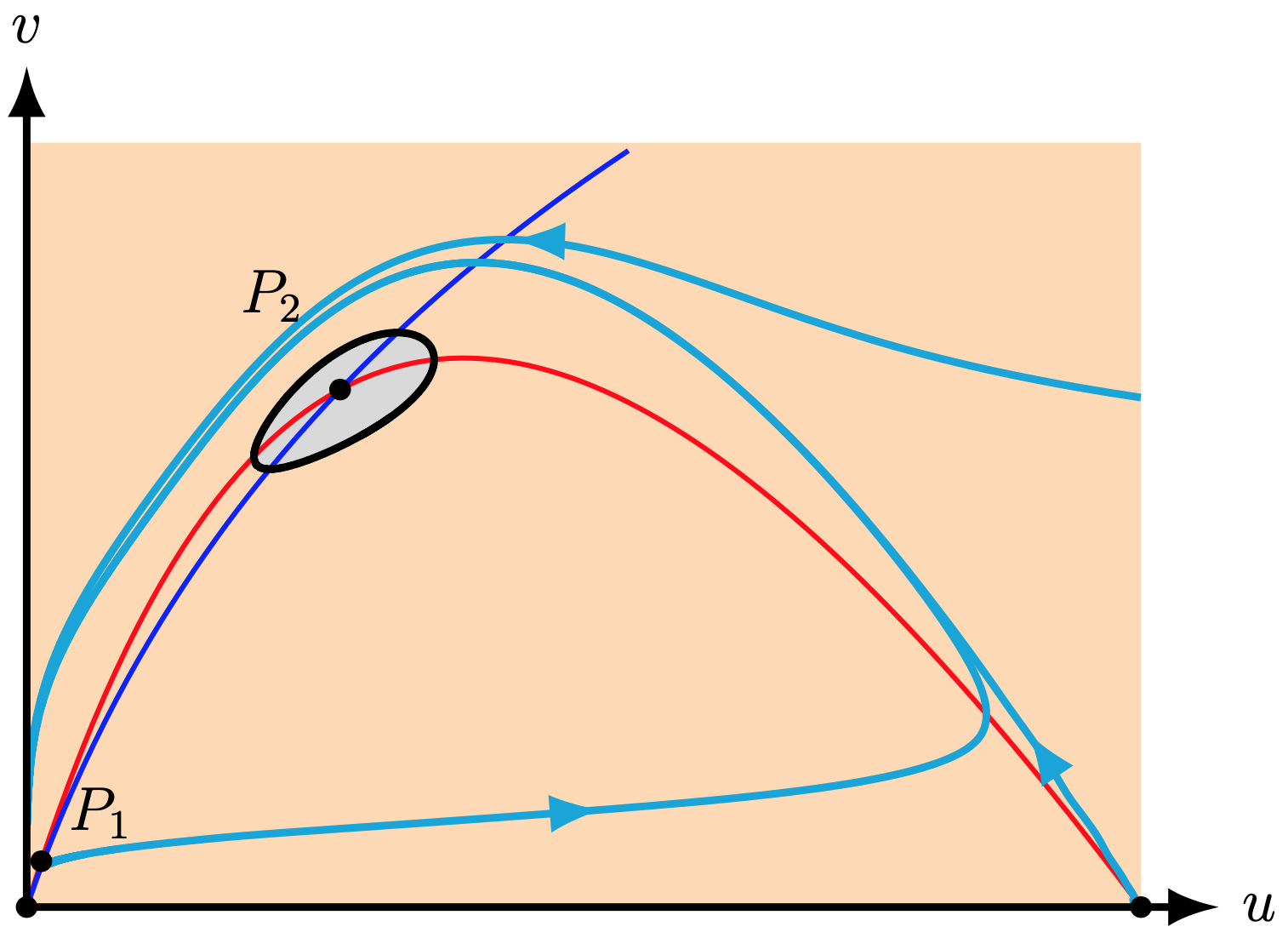}
\includegraphics[width=5.3cm]{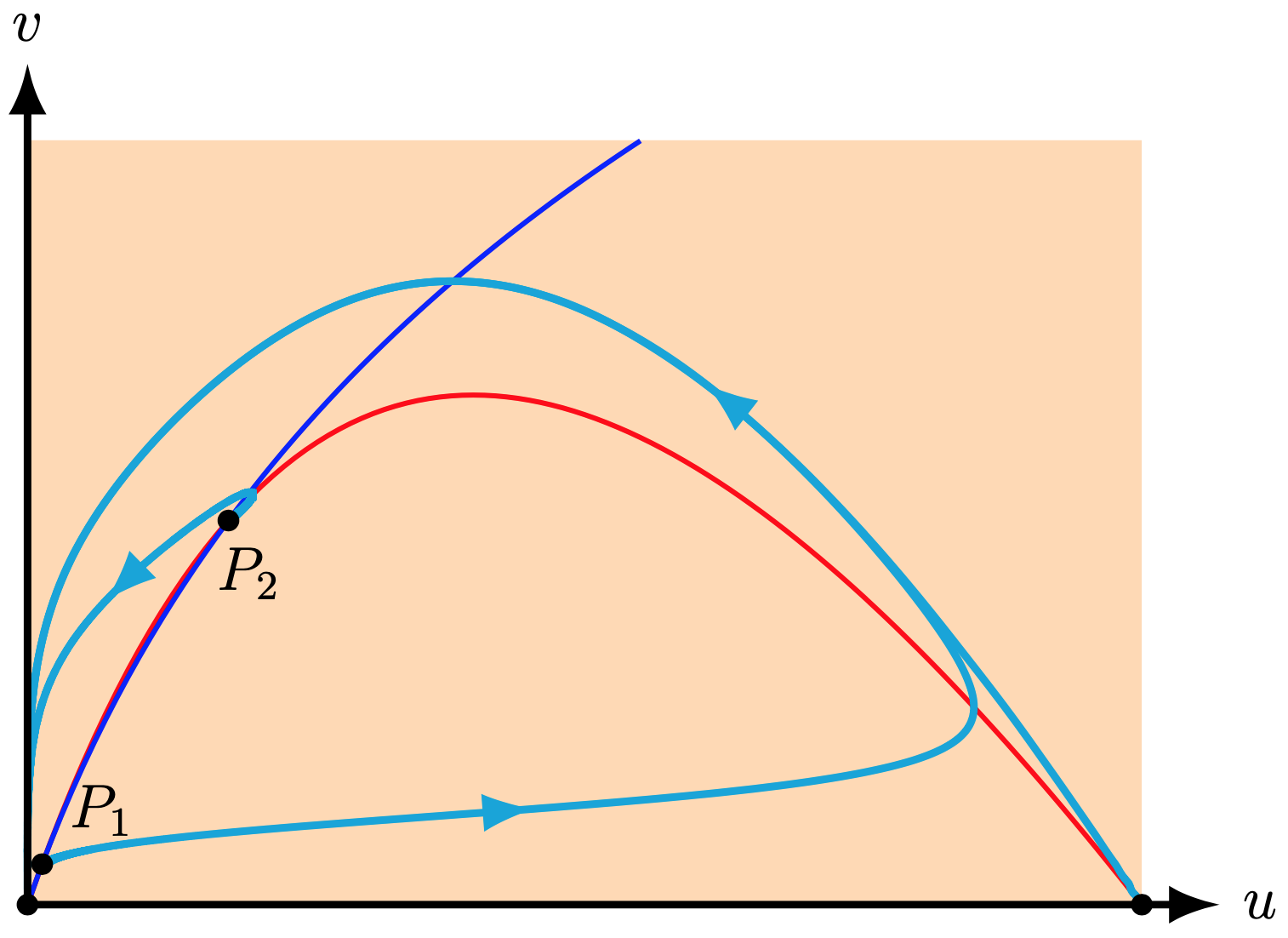}
\caption{The blue (red) curve represents the prey (predator) nullcline and the orange (grey) region represents the basin of attraction of $(0,0)$ ($P_2$). Note that $\Delta>0$ in all panels. Let $C=0.363$, $M=0.16$ and $N=0.25$, while we vary $Q$. The origin $(0,0)$ and $P_2$ are both attractors for $Q=1.695$ (top left panel) and $Q=1.69873$ (top right panel) and the equilibrium point $P_1$ is a saddle and its stable manifold (black curve) forms the boundary between the basins of attraction. For $Q=1.7$ there exist a homoclinic curve (black curve) to the saddle equilibrium point $P_1$ separating the basins of attraction of the attractors $(0,0)$ and $P_2$ (middle left panel). This homoclinic curve is broken into an unstable limit cycle upon further increasing $Q$ and the unstable limit cycle acts as a separatrix between the basins of attraction of $(0,0)$ and $P_2$ ($Q=1.705$ in the middle right panel and $Q=1.714$ in the left bottom panel). Upon further increasing $Q$, $P_2$ becomes unstable and the origin (0, 0) becomes a (global) attractor ($Q=1.8$ in the right bottom panel).}\label{Fig5} 
\end{figure}

\begin{figure}
\centering
\includegraphics[width=10.5cm]{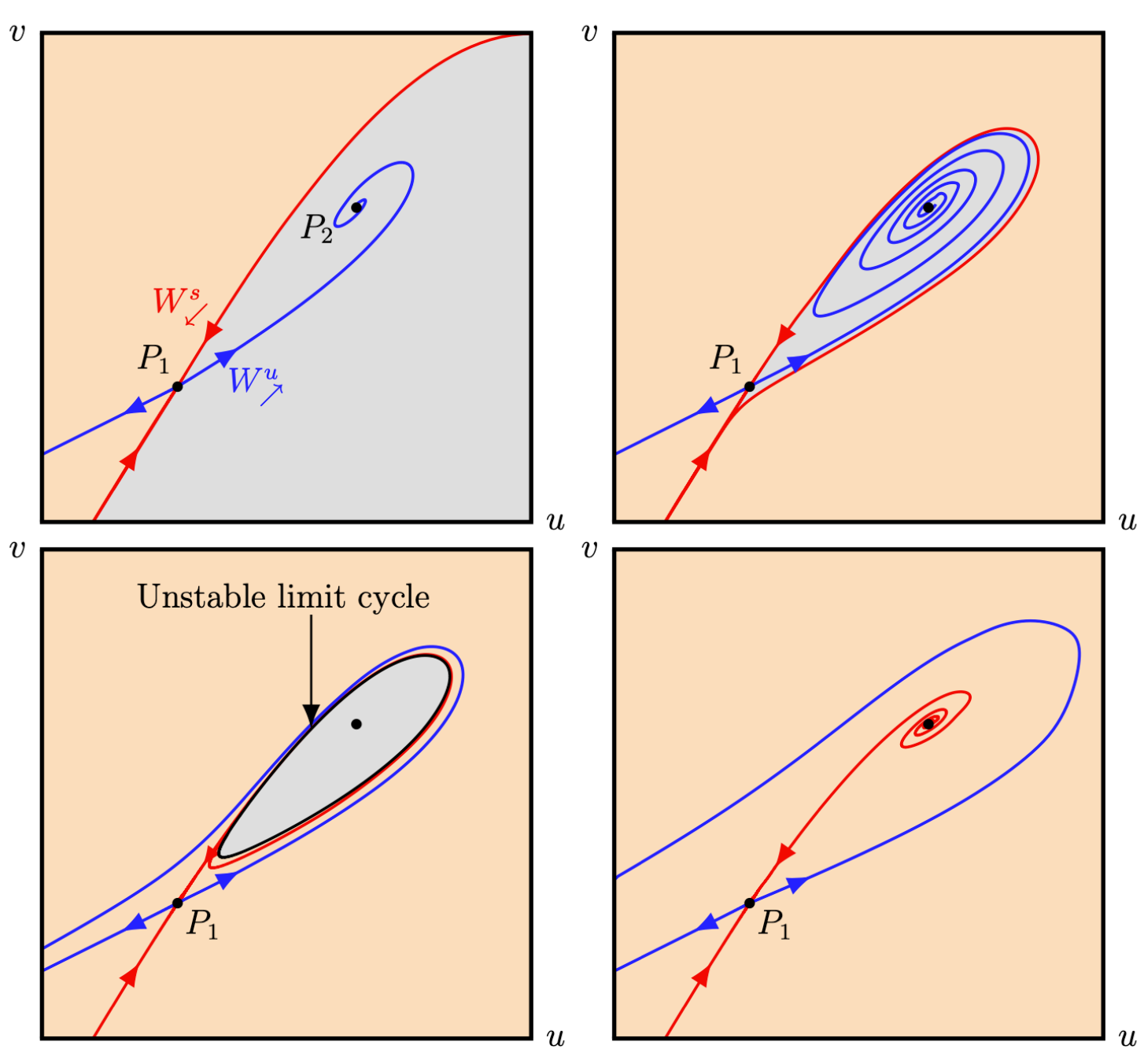}
\caption{The schematic phase planes in the neighbourhood of the equilibrium points $P_1$ and $P_2$ for increasing $Q$. The red curves represent the stable manifold of $P_1$ ($W^{s}(P_1)$) while the blue curves represent the unstable manifold of $P_1$ ($W^{u}(P_1)$. Upon increasing $Q$ we observe that the basin of attraction of the positive equilibrium point $P_2$ shrinks (light blue region), while the basin of attraction of the origin increases (orange region). We observe the birth of a homoclinic curve and its subsequent transformation into an unstable limit cycle (black curve, bottom left). An animated version of this figure is accessible on http://www.doi.org/10.6084/m9.figshare.12186708.}
\label{Fig11}
\end{figure}

Next, we discuss the stable manifold $W^s(P_1)$ of the saddle point $P_1$. This manifold often acts as a separatrix curve between the basins of attraction of the equilibrium points $(0,0)$ and $P_2$, see, for example, Figure \ref{Fig5}. 
Let $W^{u,s}_{\nearrow,\swarrow}(P_1)$ denote  
the branch of the (un)stable manifold of $P_1$ whose trajectory flows according to the direction of the arrow, see  Figure~\ref{Fig11} for an example.
Following Flores and Gonz\'alez‐-Olivares~\cite{flores} and Aguirre \textit{et al}.~\cite{aguirre2}, we get from the orientation of the nullclines and the eigenvectors of $P_1$ that $W^s_{\nearrow}(P_1)$ and $W^u_{\swarrow}(P_1)$ are connected with $(0,0)$. By continuity of the vector field in $(Q,C)$, there are parameter values for which the trajectory $W^{s}_{\swarrow}(P_1)$ intersects the boundary of $\Gamma$ (defined in Lemma~\ref{bou}). Therefore, the stable manifold $W^{s}(P_1)$ of $P_1$ act as a separatrix between the basin of attraction of the equilibrium point $P_2$ (when this is an attractor) and the equilibrium point $(0,0)$, see the top panels of Figures~\ref{Fig5} and~\ref{Fig11}. 
We observe that by changing the parameter $Q$ the stable manifold $W^s_{\swarrow}(P_1)$ moves down and the basin of attraction of $P_2$ shrinks (see light blue region of Figure~\ref{Fig5}), while the basin of attraction of $(0,0)$ increases (see orange region of Figure~\ref{Fig5}). Additionally, by the continuous dependence of the vector field on the parameter $Q$, there exist a $Q$ value for which $W^{u}_{\nearrow}(P_1)$ intersects $W^{s}_{\swarrow}(P_1)$ to form a homoclinic curve, see Figures~\ref{Fig5} and~\ref{Fig11}. Upon further decreasing $Q$, and if the equilibrium point $P_2$ is an attractor and the $\omega$-limit of the unstable manifold $W^{u}_{\nearrow}(P_1)$ is the equilibrium point $(0,0)$, the homoclinic curve breaks and there exists an unstable limit cycle which acts as a separatrix between the basins of attraction of these equilibrium points. 
Note that when the homoclinic curve breaks it generates a non-infinitesimal limit cycle (originating from a homoclinic bifurcation), which could coexist with the other limit cycle obtained via a Hopf bifurcation (infinitesimal limit cycle) when $P_2$ is a center-focus, see the upcoming Section~\ref{sec:hopf} and Figure~\ref{Fig6}. 

We summarise the above discussion in the following lemma:

\begin{lemma}
There exist conditions on the parameter values for which there is
\begin{enumerate}[label=(\roman*)]
\item a homoclinic curve determined by the stable and unstable manifold of equilibrium point $P_1=\left(u_1,v_1\right)$; and
\item a limit cycle that bifurcates from the homoclinic curve which surrounds the equilibrium point $P_2=\left(u_2,v_2\right)$.
\end{enumerate}
\end{lemma}

\subsubsection{The collapse of the positive equilibrium points}
\begin{figure}
\centering
\includegraphics[width=6.5cm]{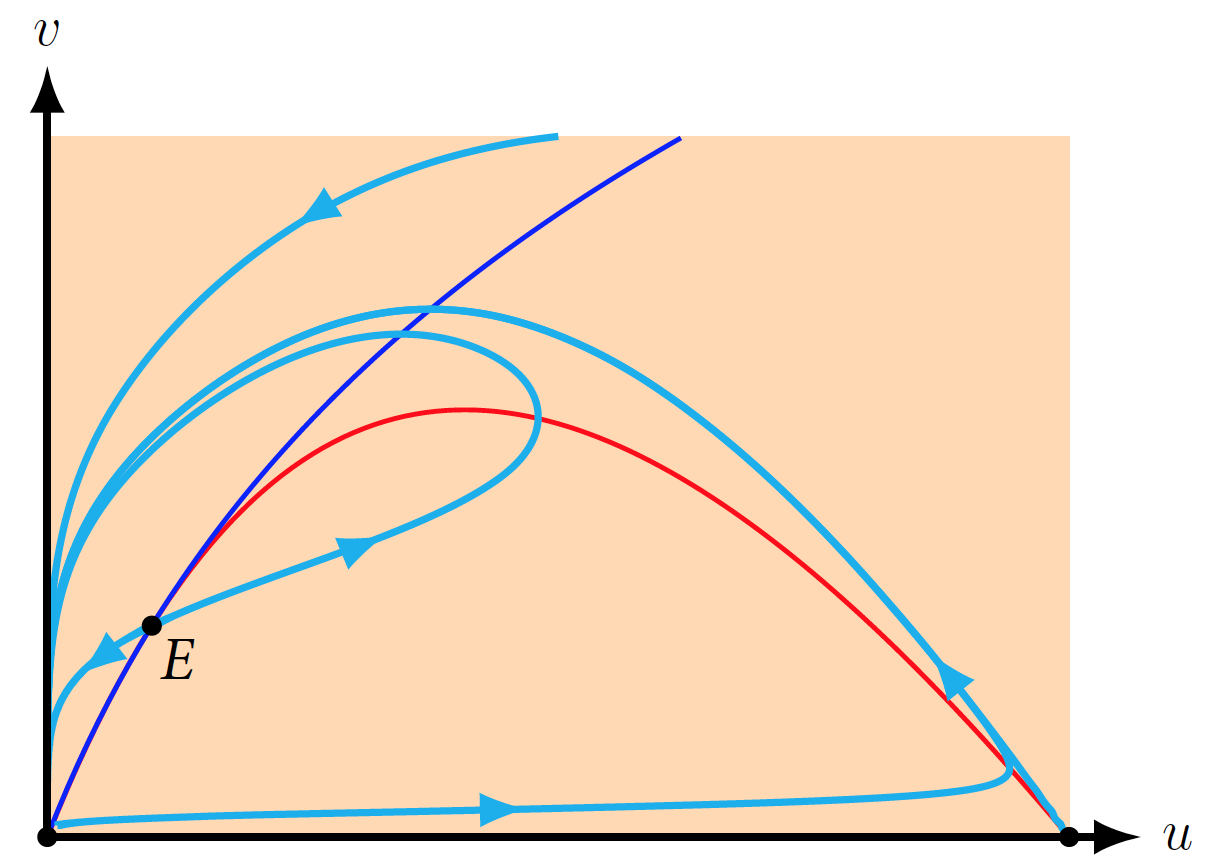}
\caption{Let $Q=1.8281$, $C=0.363$, $M=0.16$, and $N=0.25$, such that we are in case \ref{DEL0} of Section \ref{Num_eq}. The two nullclines intersect in one point in the first quadrant, i.e. $\Delta=0$~\eqref{delta}, and $P_1=P_2=E=(u_E,v_E)$. The equilibrium point is a saddle-node repeller since $C< C^*$, see Lemma~\ref{p3}. See Figure~\ref{Fig5} for the colour conventions.}
\label{Fig7}
\end{figure}

Next, we study the case~\ref{DEL0} in Section~\ref{Num_eq} where $C>M$, $\Sigma_{2}>0$, $\Delta=0$ and $N>M$ in equation~\eqref{delta}. The equilibrium points $P_1$ and $P_2$ collapse such that $u_1=u_2=u_E=(2\Sigma_2+Q(M-N))/(2(C+NQ))$ and $v_1=v_2=v_E=(2N\Sigma_2-C(M-N))/(2N(C+NQ))$. That is system~\eqref{model2} has one equilibrium point of order two in the first quadrant given by $E=(u_E,v_E)$, see Figure~\ref{Fig7}.

\begin{theorem}\label{p3}
Let the system parameter be such that $C>M$, $\Sigma_{2}>0$, $\Delta=0$ and $N>M$ (case~\ref{DEL0} of Section~\ref{Num_eq}), then the stability of equilibrium point $E=(u_E,v_E)$ is as follows:
\begin{enumerate}[label=(\roman*)]
\item $E$ is a saddle-node repeller if $0<C<\dfrac{-A+\sqrt{16MN^3(M-N)^2+A^2})}{8N^2}$\,,\\
\item $E$ is a saddle-node attractor if $C>\dfrac{-A+\sqrt{16MN^3(M-N)^2+A^2})}{8N^2}$\,,
\end{enumerate}
where $A:=-8MN^2-(1-N)(M+N)^2$.
\end{theorem}
\begin{proof}
Since $u_E+v_E=(N-M)/(2N)$ and $\Delta=0$ we have $Q=\left((M-N)^2+4CN\right)/\left(4N(C-M)\right)$. Therefore, the Jacobian matrix of~\eqref{model2} at the equilibrium point $E$ becomes   
\begin{equation}\label{jac_p1p2}
J(E)=\dfrac{M-N}{2 N^2\left(2C-M+N\right)^2}\begin{pmatrix}
 -C\left(M-N\right)^2 && \dfrac{\left(4CN+(M-N)^2\right)\left(M+N\right)^2}{4(C-M)}\\
-4CN\left(C-M\right)^2 && N(C-N)\left(4CN+(M-N)^2\right)
\end{pmatrix}.
\end{equation}
Thus, $\det{J(E)}=0$ and the trace of the Jacobian matrix~\eqref{jac_p1p2} is given by
\begin{equation}\label{qeq}
\tr\left(J(E)\right)=\dfrac{4\left(M-N\right)\left(C-M\right)}{N^2\left(2C-M+N\right)^2}\left(4N^2C^2+\beta C-\gamma\right)\,,
\end{equation}
with $\beta=N\left(M^2-6MN+N^2\right)-\left(M+N\right)^2$ and $\gamma=MN\left(M-N\right)^2$. The trace $\tr\left(J(E)\right)=0$ if and only if 
\begin{equation}\label{tr_0}
C=C^{\ast}=\dfrac{-A+\sqrt{16MN^3(M-N)^2+A^2})}{8N^2},
\end{equation}
with $A$ as given in the statement of this theorem. For $0<C<C^*$~\eqref{tr_0}, we have that $\tr\left(J(E)\right)>0$~\eqref{qeq} and the equilibrium point $E$ is a saddle-node repeller. For $C>C^*$~\eqref{tr_0}, we have that $\tr\left(J(E)\right)<0$~\eqref{qeq} and the equilibrium point $E$ is a saddle-node attractor.
\hfill\end{proof}

\subsubsection{The collapse of $P_1$ with $(0,0)$}
Next, we study the case when $C>M$, $\Sigma_2=0$ and $N>M$ such that $P_1$ merges with $(0,0)$. That is, we study case~\ref{ii} of Section~\ref{Num_eq}. Equation~\eqref{delta} has one positive solution and system~\eqref{model2} thus has one positive equilibrium point $P_2$ which is given by $P_2 =(Q(N-M)/(C+NQ)$, $C(N-M)(N(C+NQ))$. 

\begin{theorem}
If $C>M$, $\Sigma_2=0$ and $N>M$ (case~\ref{ii} of Section~\ref{Num_eq}), then the equilibrium point $P_2$ is a stable node.
\end{theorem}
\begin{proof}
Evaluating the determinant~\eqref{det} at $P_2$ gives 
\begin{equation*}
\det\left(J\left(u_2,v_2\right)\right)=\dfrac{M(2(N-M)+Q(3M^2+N(N+4M+1)))}{N}>0, 
\end{equation*} 
since $N>M$. 
The trace $\tr\left(J\left(u_2,v_2\right)\right)$~\eqref{det} gives
\begin{equation*}
\tr\left(J\left(u_2,v_2\right)\right)=-\dfrac{(N-M)^2(2NQ^2+CQ+C^2)}{N(C+NQ)^2}<0. 
\end{equation*} 
\hfill\end{proof}

Finally, system~\eqref{model2} has no positive equilibrium points in the first quadrant in the remaining cases of Section~\ref{Num_eq}, that is, in cases~\ref{DEL} and~\ref{BB}. Therefore, the equilibrium point $(0,0)$ is a global attractor in these cases since $C>M$.


\section{Bifurcation Analysis}\label{bif}
In this section we discuss some of the possible bifurcation scenarios of system~\eqref{model2}. The stability of $(1,0)$ and $P_{1}$ does not change for $\Delta>0$~\eqref{delta}, $Q>1$ and $C>M$, but whenever $\Delta=0$, $P_1$ undergoes a bifurcation as it collides with $P_2$. The stability of the equilibrium point, which is called $E$, at the bifurcation depends in a continuous fashion on the system parameters $C$ and $Q$, see Theorem~\ref{p3}.

Therefore, $C$ and $Q$ can be selected as natural candidates to act as bifurcation parameters.

\subsection{Saddle-node Bifurcation}\label{SN_S}
\begin{theorem}\label{sn}
Let the system parameters be such that $\Delta=0$~\eqref{delta}, $Q>1$, $C>M$, and 
\begin{equation}\label{sneq1}
C-3M-10N+\dfrac{4N(Q+7)(C-M)^2}{(M+N)^2}\neq0, 
\end{equation} 
then system~\eqref{model2} undergoes a saddle-node bifurcation at the equilibrium point $E$ (for changing $C$).
\end{theorem}
\begin{proof}
The proof of this theorem is based on Sotomayor's Theorem~\cite{perko}.  
For $\Delta=0$, there is only one equilibrium point $E=(u_E,v_E)$ in the first quadrant, with $u_E=(2C(1-Q)+Q(M+N))/(2(C+NQ))$ and $v_E=(CN(2Q-1)-2MNQ-CM)/(2N(C+NQ))$. 
From the proof of Theorem~\ref{p3} we know that $\det\left(J(E)\right)=0$ if $\Delta=0$. Additionally, let $U=(1,1)^T$ be the eigenvector corresponding to the eigenvalue $\lambda=0$ of the Jacobian matrix $J(E)$, and let 
\begin{align*}W=\left(-\dfrac{4N(C-M)^2}{(M+N)^2},1\right)^T\end{align*} be the eigenvector corresponding to the eigenvalue $\lambda=0$ of the transposed Jacobian matrix $J(E)^T$. 
	
If we represent~\eqref{model2} by its vector form
\begin{align} \nonumber
F(u,v;C) =\begin{pmatrix}
u(1-u)(u+v)-Quv\\ 
Cuv-v(u+v)(M+Nv)
\end{pmatrix},
\end{align}
then differentiating $F$ at $E$
with respect to the bifurcation parameter $C$ gives
\[F_C(u_E,v_E;C)=\begin{pmatrix}
0 \\ 
\dfrac{\Sigma_1\Sigma_3}{4N(C+NQ)^2}
\end{pmatrix}.\]
Therefore,
\[W \cdot F_C(u_E,v_E;C)=-\dfrac{\Sigma_1\Sigma_3}{4N(C+NQ)^2}.\]
Note that $W \cdot F_C(u_E,v_E;C)\neq0$ under certain conditions on the parameters $(C,Q,N,M)$. Next, we analyse the expression $W \cdot [D^2F_C(u_E,v_E;C)(U,U)]$ and we first compute the quadratic form associated to the Hessian matrix	
\[\begin{aligned}
D^2F(u,v;C)(V,V) =& \dfrac{\partial^2F(u,v;Q)}{\partial u^2}v_1v_1+\dfrac{\partial^2F(u,v;Q)}{\partial u\partial v}v_1v_2+\dfrac{\partial^2F(u,v;Q)}{\partial v\partial u}v_2v_1+\dfrac{\partial^2F(u,v;Q)}{\partial v^2}v_2v_2\,.
\end{aligned}\]
At the equilibrium point $E$ and $U$, this becomes
\[\begin{aligned}	
D^2F(u_E,v_E;C)(U,U)& = \begin{pmatrix}
-(Q+7)\\ 
C-3M-10N
\end{pmatrix}\,. 
\end{aligned}\]
Therefore, 
\[\begin{aligned}
W \cdot [D^2F(u_E,v_E;C)(U,U)]=C-3M-10N+\dfrac{4N(Q+7)(C-M)^2}{(M+N)^2} \,,
\end{aligned}\]
and $W \cdot [D^2F(u_E,v_E;C)(U,U)]\neq0$ by our assumptions on the parameters $(C,Q,N,M)$~\eqref{sneq1}. Therefore, by Sotomayor's Theorem~\cite{perko} it now follows that system~\eqref{model2} has a saddle-node bifurcation at the equilibrium point $E$ if the conditions on the parameters are met.
\hfill\end{proof}

\subsection{Hopf Bifurcation}\label{sec:hopf}
Let $(U,V)$ be the coordinates of an equilibrium of~\eqref{model2} in the first quadrant. Then the set of equations $W(U,V)=0$ and $R(U,V)=0$
define implicitly a locally invertible transformation given by
\[\begin{aligned}
\Psi&:\Lambda \longrightarrow\mathbb{R}^4_+,\\
(C,M,U,V)\mapsto(C,M,N,Q)&:=\left(C,M,\dfrac{CU-M(U+V)}{V(U+V)},\dfrac{(1-U)(U+V)}{V}\right),
\end{aligned}\]
in $\Lambda:=\{(C,M,U,V)\in\mathbb{R}^4_+:\,CU-M(U+V)>0,\,\, U-1<0\}.$
Considering the time rescaling $\tau\mapsto\tau/(V(U+V))$, system~\eqref{model2} in parameter space $\Lambda$ has the form:
\begin{equation}\label{hopf1}
\begin{aligned}
\dfrac{dx}{d\tau} &= x(U+V)\left(Vx-Vx^2-Uy+U^2y+UVy-Vxy\right), \\
\dfrac{dy}{d\tau} &= y\left(CV^2x-M(U+V)(V-y)(x+y)+CU(Vx-y(x+y)\right) \,, 
\end{aligned}
\end{equation}
where, for convenience, we again use the notation $(x,y)$ to name the state variables. System~\eqref{hopf1} is $C^{\infty}$-equivalent to~\eqref{model2} in parameter space $\Lambda$. Moreover, the (positive) equilibrium coordinates appear now as the explicit parameters $(U,V)$. 
In what follows, we will derive conditions such that \eqref{hopf1} undergoes a generic Hopf bifurcation at $(U,V)$, see~\cite{guckenheimer,kuznetsov} for more details.

The Jacobian matrix of~\eqref{hopf1} at $(U,V)$ is
\begin{equation}J(U,V)=\begin{pmatrix} \label{jacUV}
-UV (U+V) (-1 + 2 U + V) & (U-1)(U+V) U^2 \\
C V^3 & V \left(M (U+V)^2 - C U (U+ 2 V)\right) \end{pmatrix}.
\end{equation}
The trace and determinant of $J(U,V)$ are given by 
\begin{equation*}
\rm{tr}\left(J(U,V)\right)=VT, \hspace{5mm} {\rm and}\hspace{5mm} \det\left(J(U,V)\right)=UV^2(U+V)^2D,
\end{equation*}
respectively, where
\begin{equation}\label{det2}
\begin{aligned}
D:=&D(C,M,U,V)=CU(-1+2U+2)+M(U-2U^2+V-3UV-V^2)~\text{and}\\
T:=&T(C,M,U,V)=\left(M(U+V)^2-CU(U+2V)-U(U+V)(-1+2U+V)\right).
\end{aligned}
\end{equation}
Whenever $T=0$ and $D>0$, the eigenvalues of $J(U,V)$ are purely imaginary and non-trivial. Moreover, since we have
\begin{equation}\label{traceM}
\dfrac{\partial T}{\partial M}(C,M,U,V)=(U+V)^2>0,
\end{equation}
the Hopf bifurcation in~\eqref{hopf1} is generically unfolded by parameter $M$~\cite{kuznetsov}. In particular, it follows from~\eqref{traceM} that equation $T(C,M,U,V)=0$ implicitly defines the function
\begin{equation}\label{MT}
M(C,U,V)=\dfrac{CU(U+2V)}{(U+V)^2}+\dfrac{U(-1+2U+V)}{U+V}.
\end{equation}

We now calculate the first Lyapunov quantity~\cite{guckenheimer,kuznetsov} in order to determine genericity conditions. We follow
the derivation in~\cite{guckenheimer} and move the equilibrium point $(U,V)$ of the system~\eqref{hopf1} to the origin via the translation $x\mapsto x+U$, $y\mapsto y+V$ to obtain the equivalent system
\begin{equation}\label{hopf2}
\begin{aligned}
\dfrac{dx}{d\tau} &=	 (U + V) (U + x) \left(-V^2 x + (-1 + U) U y - V x (-1 + 2 U + x + y)\right), \\
\dfrac{dy}{d\tau} &=	 (V + y) \left(M (U + V) (U + V + x + y)y + C (V^2 x - 2 U V y - U y (U + x + y))\right) \,.
\end{aligned}
\end{equation}
In particular, the Jacobian matrix of~\eqref{hopf2} at the equilibrium point $(0,0)$ coincides with $J(U,V)$ in \eqref{jacUV}. Substitution of~\eqref{MT} into $J(U,V)$ gives
\begin{equation*}
J_H(U,V)=\begin{pmatrix} 
-UV (U+V) (-1 + 2 U + V) & (U-1)(U+V) U^2 \\
C V^3 & UV \left(2 U^2 + (-1 + V) V + U (-1 + 3 V))\right) \end{pmatrix},
\end{equation*}
with $\tr\left(J_H\left(U,V\right)\right)\equiv0$ and $\det\left(J_H(U,V)\right)=UV^2(U+V)^2 D_H$~\eqref{det2} where
\begin{equation*} 
D_H=D|_{T=0}=-U - 4 U^3 + C V - (C-6) UV - (V-1)^2 V - 5 UV^2 - U^2 (8 V-4).
\end{equation*}
If $D_H>0$, then $\mathbf{v}_1=\left(\begin{array}{c} \dfrac{U (U - 2 U^2 + V - 3 UV - V^2)}{CV^2} \\ 1 \end{array} \right)$ and $\mathbf{v}_2=\left(\begin{array}{c} -w \\  0 \end{array} \right)$ are the generalised eigenvectors of $J_H(U,V)$, where 
\begin{equation}\label{valorw}
w=V(U+V)\sqrt{UD_H}.
\end{equation}
The change of coordinates $\left(\begin{array}{c} x \\  y \end{array} \right)\mapsto [\mathbf{v}_{1} \, \mathbf{v}_{2}] \left(\begin{array}{c} x \\  y \end{array} \right) $ allows us to express system (\ref{hopf2}) with $T(C,M,U,V)=0$ in the form
\begin{equation}\label{sisparen}
\left(\begin{array}{c}
\dfrac{dx}{d\tau}\\ 
\\
\dfrac{dy}{d\tau}
\end{array} \right)=\left( \begin{array}{cc}
0 & -w\\ 
\\
w & 0
\end{array} \right)\left(\begin{array}{c}
x \\ 
\\
y
\end{array} \right)+\left(\begin{array}{c}
P(x,y) \\ 
\\
Q(x,y)
\end{array} \right),
\end{equation}
where
\begin{equation}\label{valorPyQ1}
\begin{aligned}
P\big(x,y\big)=&-\dfrac{x^3}{C V^2}\big(C U - M U - M V\big) \big(G U + C V^2\big)+\dfrac{x}{CV^2}\big(CGUV^3-C^2U^2V^3+C M U^2 V^3 - 2 C^2 U V^4\\
& + 2 C M U V^4+ C M V^5+ C U w y - M U w y- C V w y - M V w y\big)-\dfrac{x^2 }{C V^3}\big(C G U^2 V^2- G M U^2 V^2 \\
&- C G U V^3- G M U V^3 + C^2 U^2 V^3 -C M U^2 V^3 + 3 C^2 U V^4- 3 C M U V^4- 2 C M V^5- C U w y\\
& + M U w y + M V w y\big),
\end{aligned}
\end{equation}
\begin{equation}\label{valorPyQ3}
\begin{aligned}
Q\big(x,y\big)=&\dfrac{G U x^3}{C^2 V^2 w}\big(G U + C V^2\big) \big(G U^2 - C^2 U V + G U V + C M U V + C M V^2\big)+\dfrac{G U x^2}{C^2 V^3 w}\big(-C G U^2 V^3 \\
&- C^2 G U^2 V^3 + C G M U^2 V^3 + 3 C G U^3 V^3 - C G U V^4 + C^2 G U V^4 + C G M U V^4+ C^2 U^2 V^4\\
& - C^3 U^2 V^4+ 4 C G U^2 V^4 + C^2 M U^2 V^4 - C^2 U^3 V^4 + C^2 U V^5 - 3 C^3 U V^5+ C G U V^5\\
&+ 3 C^2 M U V^5 + 2 C^2 M V^6 + C^2 U V^6 - 3 G U^2 W y + w y\big(C^2 U V- 3 G U V - C M U V- C M V^2\\
&- 2 C U V^2 - 2 C V^3\big)\big)-\dfrac{wy^2}{C^2V^5}\big(C U V^3 - 3 C U^2 V^3 + C V^4 - 4 C U V^4 - C V^5 + U w y + V w y\big)\\
&+\dfrac{xy}{C^2 V^4 }\big(2 C G U^2 V^3 + C^2 G U^2 V^3- C G M U^2 V^3 - 6 C G U^3 V^3  + 2 C G U V^4  - C^2 G U V^4 \\
&- C G M U V^4 - C^2 U^2 V^4 - 8 C G U^2 V^4 + C^2 U^3 V^4- C^2 U V^5 - 2 C G U V^5- C^2 U V^6 + 3 G U^2 w y \\
&+ 3 G U V w y+ C U V^2 w y + C V^3 w y\big),\\
\end{aligned}
\end{equation}
and $G=G\big(U,V\big)=U - 2 U^2 + V - 3 U V - V^2.$
 
System (\ref{sisparen}) and equations \eqref{valorw},~\eqref{valorPyQ1} and~\eqref{valorPyQ3}  allow us to use the derivation in~\cite{guckenheimer} for the direct calculation of the first Lyapunov quantity $L_1$. In this way we obtain the following expression:
\begin{equation*}
L_1=\dfrac{ U^3}{8 C (U + V) w^2}\,l_1,
\end{equation*}
where
\begin{equation}\label{valorPyQ2}
\begin{aligned}
 l_{1}=&-C^3 V^4 \big(-U + U^2 - 2 U V - 2 V^2\big) -U \big(1 - U\big) \big(U + V\big)^3 \big(-2 U^3 + 4 U^4 + V- 9 U V + 20 U^2 V\\
& - 10 U^3 V- 6 V^2 + 28 U V^2 - 28 U^2 V^2 + 9 V^3 - 19 U V^3- 4 V^4\big)- C^2 V^3 (3 U^2 - 12 U^3 + 12 U^4 + 3 U V\\
& - 14 U^2 V + 16 U^3 V + 2 U V^2 + 4 V^3 - 13 U V^3 - 4 V^4) + C V (U + V)^2 \big(-2 U^4 + 2 U^5 + 3 U V - 21 U^2 V \\
&- 5 U^2 V^2+ 41 U^3 V - 27 U^4 V - 10 U V^2 + 38 U^2 V^2 - 38 U^3 V^2+ 2 V^3 - 8 U^2 V^3 - 4 V^4 + 6 U V^4 + 2 V^5\big).\\
\end{aligned}
\end{equation}

Thus we have obtained the following result.

\begin{theorem}\label{hopfp3}
Let $(C,M,U,V)\in \Lambda$ be such that $T(C,M,U,V)=0$, $D_H>0$ and $l_1\neq0$~\eqref{valorPyQ2}.
Then (\ref{hopf1}) undergoes a codimension-one Hopf bifurcation at the equilibrium point $(U,V)$. In particular, if $l_1<0$ (resp. $l_1>0$), the Hopf bifurcation is supercritical (respectively subcritical), and a stable (respespectively unstable) limit cycle bifurcates from $(U,V)$ under suitable parameter variation.
\end{theorem}

\begin{figure}
\centering
\includegraphics[width=15cm]{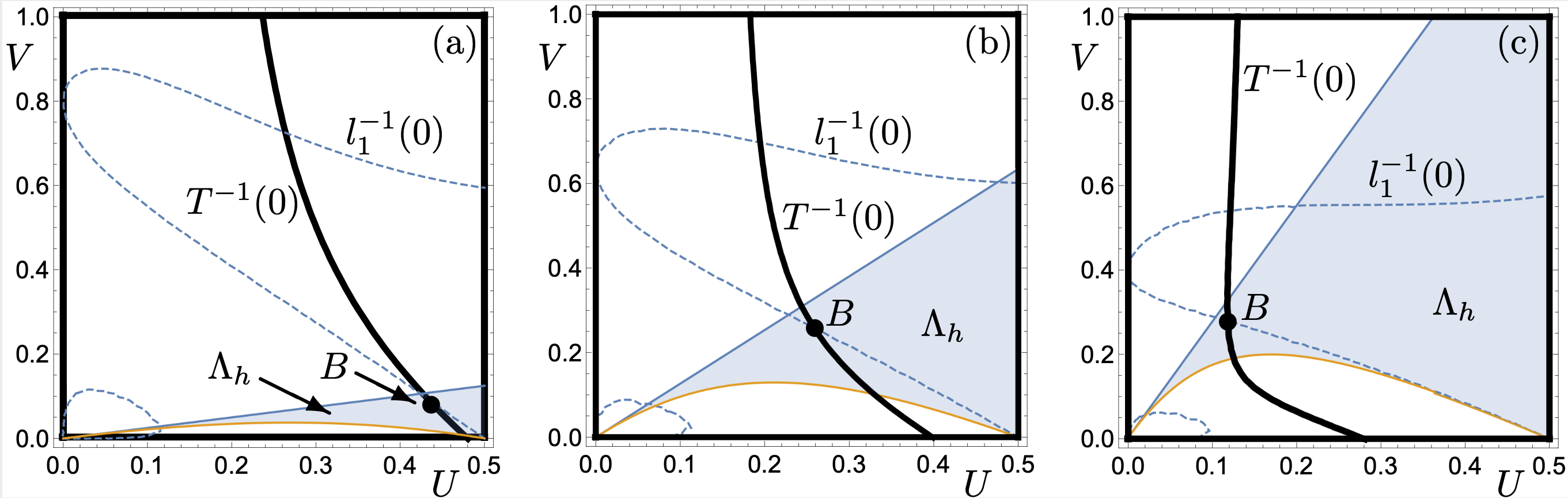}
\caption{Bifurcation sets in the $(U,V)$-parameter plane. The black curve defined by $T(C,M,U,V)=0$ corresponds to a Hopf bifurcation in the shaded region labelled as $\Lambda_h$. Parameter values are $C=0.2$ in panel (a), $C=0.363$ in panel (b), and $C=0.6$ in panel (c), while $M=0.16$ is fixed. }
\label{Fig8}
\end{figure}

\begin{figure}
\centering
\includegraphics[width=14cm]{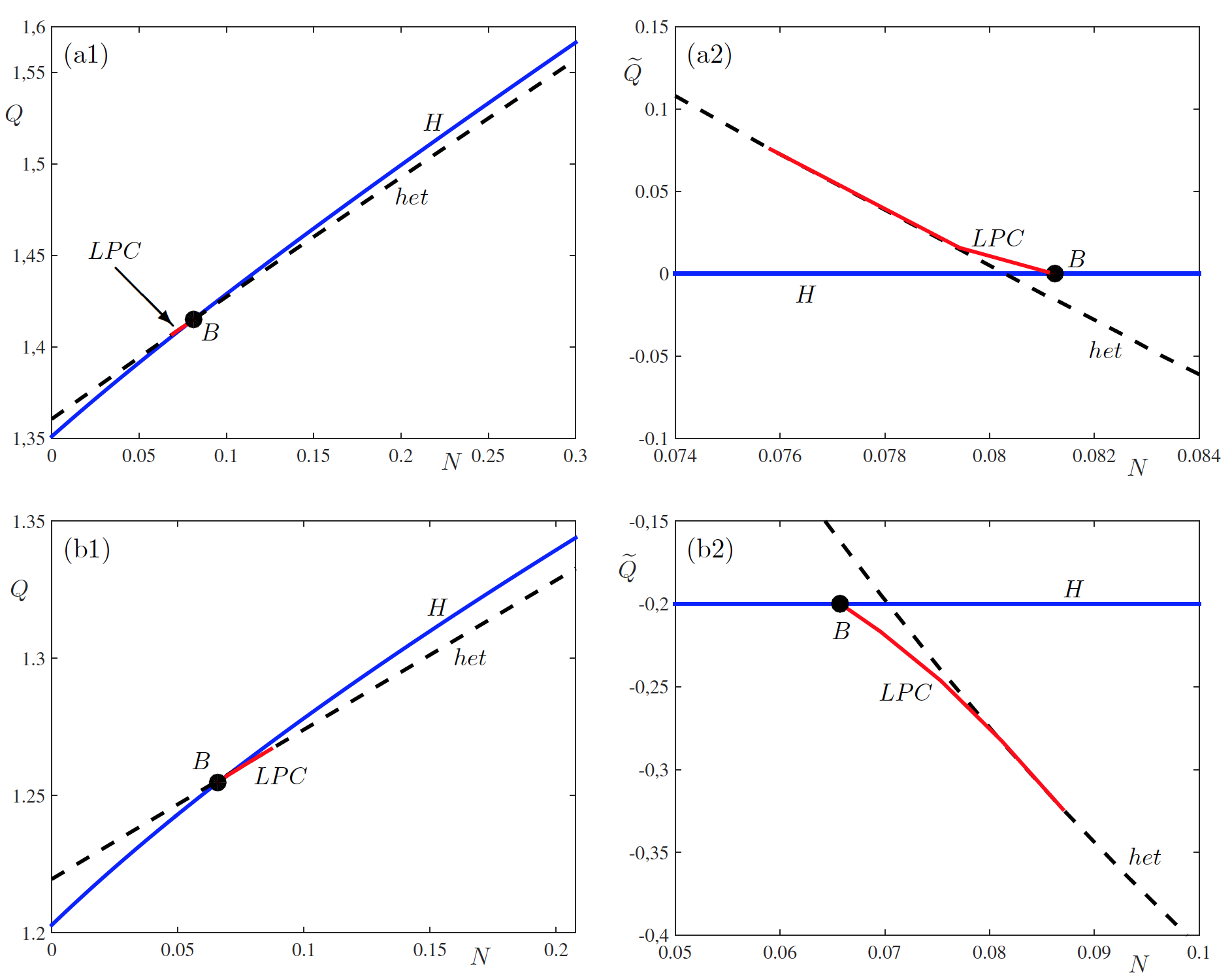}
\caption{Bifurcation diagrams in the $(N,Q)$-plane (left panels) and $(N,\widetilde{Q})$-plane (right panels) in a neighbourhood of the Bautin point $B$ when $M=0.16$ is fixed and $C=0.363$ (top panels) and $C=0.6$ (bottom panels). $LPC$ represent a curve of saddle-node (or limit point) bifurcations of limit cycles, $het$ represent a heteroclinic bifurcation, $H$ represent a Hopf bifurcation. }
\label{Fig12} 
\end{figure}

\begin{figure}
\centering
\includegraphics[width=8cm]{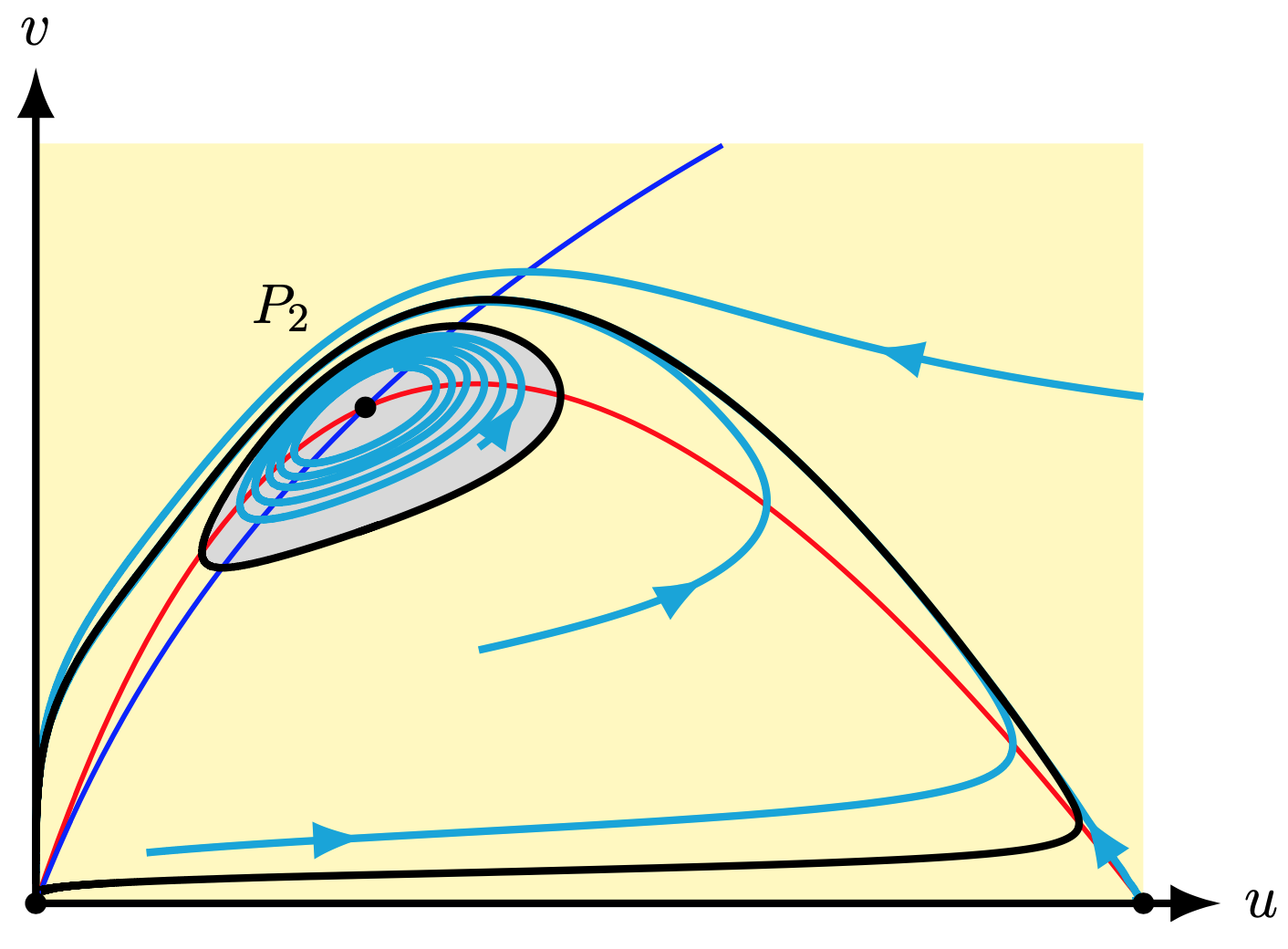}
\caption{Let $Q=1.77$, $C=0.363$, $M=0.17$ and $N=0.25$ such that we are in case~\ref{i} of section \ref{Num_eq}. The origin $(0,0)$ is an unstable node and the equilibrium point $P_2$ is stable surrounded by two limit cycles. The unstable limit cycle act as a separatrix between the basins of attraction of the stable limit cycle and $P_2$. The yellow (grey) region represent the basin of attraction of the stable limit cycle ($P_2$). See Figure~\ref{Fig5} for the colour conventions.}
\label{Fig6} 
\end{figure}

Figure~\ref{Fig8} shows the relevant bifurcation sets in the $(U,V)$-plane for increasing values of $C$. The black curve defined by $T(C,M,U,V)=0$ corresponds to a Hopf bifurcation only inside the shaded region labelled as $\Lambda_h=\Lambda\cap\{D_H>0\}$.
The Hopf bifurcation curve is separated into two segments: the upper one corresponds to $l_1<0$ (supercritical), and the lower one to $l_1>0$ (subcritical). The division occurs at the point $B=L_1^{-1}(0)\cap T^{-1}(0)\cap \Lambda_H$ --also known as Bautin bifurcation-- where the Hopf bifurcation curve intersects the level set $\{L_1=l_1=0\}$; at such point the Hopf bifurcation of (\ref{hopf1}) at $(U,V)$ is degenerate. The actual codimension of this singularity -- and the stability of further limit cycles that bifurcate -- is determined by the sign of the so-called second Lyapunov quantity $L_2$~\cite{kuznetsov}. In order to check the sign of $L_2$ for representative parameter values, we performed a bifurcation analysis with Auto-07p~\cite{auto}. Figure~\ref{Fig12} shows the corresponding bifurcation diagrams in the $(N,Q)$-plane in a neighbourhood of the Bautin point $B$ when $C=0.363$ and $C=0.6$ in rows (a) and (b), respectively, and $M=0.16$ remains fixed. The right column shows the same curves near $B$ after the change of parameter $\widetilde{Q}=100(Q-H(N))$, where $Q=H(N)$ defines (locally) the curve of Hopf bifurcation in the $(N,Q)$-plane. The advantage of this rescaling is that the curves in the $(N,\widetilde Q)$-plane are now well separated. In Figure~\ref{Fig12}, a curve of saddle-node (or limit point) bifurcations of limit cycles ($LPC$) emanates from $B$ and eventually moves towards a curve of heteroclinic bifurcation (labeled $het$). At the curve $het$, a heteroclinic cycle is formed which joins $(0,0)$ with $(1,0)$ along the horizontal axis, and then goes back to the origin in the interior of the phase space. Moreover, for parameter values in the small region bounded by the curves $LPC$ and $het$, we have two concentric limit cycles surrounding the point $P_2$ in phase space. The continuation procedure confirms that for $C=0.363$ we have $L_2>0$ and the largest limit cycle is unstable (similar to Figure~\ref{Fig5}); and for $C=0.6$, we have $L_2<0$, so the largest limit cycle is stable (similar to Figure~\ref{Fig6}). The theory of generalised Hopf bifurcations is only local, meaning that it is valid for parameter values in a neighbourhood of the special Bautin point $B$~\cite{guckenheimer,kuznetsov}. Hence the two phase portraits in Figures~\ref{Fig5} and~\ref{Fig6} do not contradict the theory nor one another: they can be found well away from each other in parameter space. We conjecture that there is a higher codimension bifurcation point organising the overall dynamics in our model (which is yet to be found), and the two phase portraits are just two possible ``perturbations" away from that point.

\subsection{Bogdanov-Takens bifurcation}\label{sec:bt}

We now give conditions such that our model undergoes a Bogdanov-Takens bifurcation under suitable parameter variation at the equilibrium point $E$. We refer to~\cite{kuznetsov} and the references therein for the derivation of the genericity and transversality conditions that need to be verified during this proof.

For the sake of clarity, it is convenient to state the dependence of the vector field~\eqref{model2} on parameters $C$ and $Q$ explicitly. Hence, throughout this section we denote
\[\begin{aligned}
& X: \mathbb{R}^4_+\longrightarrow\mathbb{R}^2_+, \\
& X(u,v;C,Q) = \bigg( u\left(1-u\right)\left(u+v\right)-Quv, Cuv-v\left(u+v\right)\left(M+Nv y\right) \bigg),
\end{aligned}\]
where we use notation $(u,v)$ for the state variables.
Also, let us denote the Jacobian matrix of $X$ with respect to the variables $(u,v)$ as
\begin{equation*}
\dfrac{\partial X}{\partial(u,v)}(u,v;C,Q).
\end{equation*}
Note that this notation for a Jacobian is standard in bifurcation theory as it explicitly highlights the dependence on all the dependent and independent variables \cite{kuznetsov}.

{\em Step 1.} We verify that the singularity has a double zero eigenvalue with geometric multiplicity one. 

From the proof of Theorem~\ref{p3}, conditions $\Delta=0$ in~\eqref{delta} and $tr(J(E))=0$ in~\eqref{qeq} determine a bifurcation point $(C,Q)=(C^*,Q^*)$ 
implicitly in the form:
\begin{eqnarray*}
C^{\ast}&=&\dfrac{-A+\sqrt{16MN^3(M-N)^2+A^2}}{8N^2} ;\\
Q^{\ast}&=&\dfrac{M^2 + 4 C^* N - 2 M N + N^2}{4 N(C^* - M )},\\\nonumber
\end{eqnarray*}
where $A=-8MN^2-(1-N)(M+N)^2.$ Moreover, at $(C,Q)=(C^*,Q^*)$, the equilibrium point $E=(u_E,v_E)$ has a Jacobian matrix given by
\begin{equation}\label{eq:jac}
\dfrac{\partial {X}}{\partial(u,v)}(u_E,v_E;C^*,Q^*)=
 \left(
\begin{aligned}
 -\dfrac{(C^* (M - N) (M + N)^2}{2 N^2 (2 C^* - M + N)^2)} & \dfrac{J_{12}(M - N) (M + N)^2}{8 (C^* - M) N^2 (2 C^* - 
  M + N)^2} \\
 \dfrac{2 C^* (C^* - M)^2 (-M + N)}{N (2 C^* - M + N)^2} & \dfrac{J_{12}(C^* - M) (M - N)}{2 N (2 C^* - M + N)^2} \\
\end{aligned}
\right)
\end{equation}
with $J_{12}=M^2 - 2 M N + N (4 C^* + N)$ and a double zero eigenvalue. However, note that~\eqref{eq:jac} is not the null matrix. The corresponding generalised eigenvectors of~\eqref{eq:jac} are given by 
\begin{equation} \label{eq:eigenvectors}
	\mathbf{v}_1=(z_1,1)^T, \hspace{5mm} \mathrm{and} \hspace{5mm}\mathbf{v}_2=\left(1,z_2\right)^T,
\end{equation}
where
\begin{eqnarray} \label{eq:v1}
	z_1&=&\dfrac{4N^3}{w_3} \left(w_1 + w_2\right),\\\label{eq:v2}
	z_2&=&\dfrac{2 N (2 C^* - M + N)^2 + 
 4 C^* (C^* - M)^2 (M - N)}{(C^* - M) (M - N) (M^2 - 2 M N + N (4 C^* + N))},\\\nonumber
\end{eqnarray}
and
\begin{align*} \nonumber
w_1=&\big(1 + N\big) \big(M + N\big)^2,\\
w_2=& \sqrt{16 M N^3 \big(M - N\big)^2 + A^2}.\\
w_3=& M^3\big (M \big(-1 + N\big)^2 + 4 N \big(1 + N\big)\big) + 2 M N \big(1 + N\big) \big(2 N^2 + w_2\big)- \big( N-1\big) N^2\big (N^2 - N^3 + w_2\big) \\
&-M^2 \big(2 N^2 \big(-3 - 6 N + N^2\big)- w_2+ N w_2\big),
\end{align*}
It follows that~\eqref{eq:jac} is nilpotent and that the double zero eigenvalue has geometric multiplicity one.

{\em Step 2.}
The next goal is to state the transversality condition of a Bogdanov-Takens bifurcation, namely, that the map 
\begin{equation*}
\Psi:\mathbb{R}^4\rightarrow\mathbb{R}^4,\,(u,v,C,Q)\mapsto(X,\mathcal{T},\mathcal{D})
\end{equation*}
is regular at $(u,v,C,Q)=(u_E,v_E,C^{\ast},Q^{\ast})$, where $\mathcal{T}$ and $\mathcal{D}$ are the trace and determinant of the Jacobian matrix $\partial X/\partial(u,v)\left(u,v;C,Q\right)$, respectively. 

After some calculations, the determinant of the $4\times4$ Jacobian matrix $D\Psi$ of the map $\Psi$ can be written as 
\begin{equation*}\det D\Psi(u,v,C,Q)=-uvF(u,v,C,Q)\end{equation*}
with:
\begin{equation*}
\begin{aligned}
	F(u,v,C,Q)=&C^2 u^3 (-2 + 6 u + v)+ M^2 u (6 u^3 - 2 (Q-1) v^2 + u^2 (7 v-2))+ N v^2 (24 (1 + N) u^4+u^2 (6 ( Q-1) \\
	&+ N (7+ 11 Q - 48 v)) v+ u^3 (3 Q + 2 N (Q - 12 v-7) + 24 v-17)-3 (Q-1) v^2 (2 Q - 3 N v-2) \\
	&-uv( Q-1)(Q + 9 v + 6 N v-7))+C u ((9 + 6 N) u^4+ 3 N ( Q-1) v^3+u v (2 Q + 3 N (v-2) v-2) \\
	&+ u^2 v (1 - 3 Q + v+6 N (1 + v))-2 M (6 u^3 - ( Q-1) v^2 + u^2 ( 4 v-2)) - u^3 (6 - 6 v + N (4 + 15 v)))\\
	&+M (-3 (3 + 2 N) u^5 - (2 + 5 N) ( Q-1) u v^3+u^3 v ( 3 Q + 2 N (v-3) + 5 v-1)-2 (Q-1) v^3 ( Q \\
	&- 3 N v-1) + u^4 (6 - 6 v + N (4 + 15 v))+ u^2 v (2 + 2 (N-3) v - 25 N v^2 + Q ((6 + 4 N) v-2))).
\end{aligned}
\end{equation*} 
In particular, straightforward substitution and algebraic simplification leads to $F(u_E,v_E,C^{\ast},Q^{\ast})=G_1$, where
\begin{equation}\label{eq:G1}
\begin{aligned} 
   	G_1=& \bigg(2 (M - N)^5 (M + N)^4 \big(M^5 (-1 + N)^3 (-1 + 3 N)+M^4 N (5 - 4 N- 18 N^2+4N^3+13 N^4) +N^3 (1 \\
	&+ N - 9 N^2 + 7 N^3) (N^2 - N^3 + Z +M N^2 (1 + N)^2(5 N^2+ 4 N^3 - 17 N^4 + 3 Z + 3 N Z)\\
	&-M^2 N (-44 N^4 + 20 N^5 + 6 N^6 - 3 Z - 3 N Z+N^3 (-36 + 7 Z) + N^2 (-10 + 7 Z))-M^3 ( 5 N Z\\
	&-8 N^4 - 40 N^5 - 14 N^6- Z +3 N^3 ( Z-8)-N^2 (10 + 7 Z)))\big)\bigg)\\
      	&/\bigg(N (N^2 -M^2 (N-1) + 3 N^3 + 2 M N (1 + N) + Z)^5\bigg),
\end{aligned}
\end{equation}
and
$Z=(M + N)\sqrt{M^2 (N-1)^2 + (N-1)^2 N^2 + 2 M N (1 + 6 N + N^2)}.$

It follows that if 
\begin{equation}\label{eq:bt1}
G_1\neq0,
\end{equation}
 then we have
$
\det D\Psi(u_E,v_E,C^{\ast},Q^{\ast})\neq0,
$
which ensures that the map $\Psi$ is regular at $(u,v,C,Q)=(u_E,v_E,C^{\ast},Q^{\ast})$.

{\em Step 3.} 
We now construct a change of coordinates ---and give sufficient conditions to do so--- which transforms $X(u,v;C,Q)$ into a normal form of the Bogdanov-Takens bifurcation; we refer to~\cite{kuznetsov} again.

Let us first define the following auxiliary expressions:
\begin{eqnarray}\label{eq:G2}
G_2&=&Q^* -1 - M + 2 u_E - N u_E - z_2 + C^* z_1-M z_1 + 3 u_E z_1- 3 N v_E+ z_1 v_E - 2 N z_1 v_E),\\ \label{eq:G3}
 G_3&=&\big( 2 u_E -1 + Q^* + 2 N u_E - 2 z_1 - C^* z_1 + 6 u_E z_1 + M (2 + z_1)+6 N v_E + 2 v_E z_1 + 2 N v_Ez_1\big)\\ \label{eq:G4}
G_4&=&z_1 z_2-1,
\end{eqnarray}
where $z_1$ and $z_2$ are given by~\eqref{eq:v1} and~\eqref{eq:v2}, respectively.

Let us now move the equilibrium point $(u_E,v_E)$ of~\eqref{model2} to the origin via the translation $u\mapsto u+u_E$, $v \mapsto v+v_E$ to obtain the equivalent system
\begin{equation}\label{eq:Y}
\begin{aligned}
\dfrac{du}{d\tau} &=(u_E + u) \big( (1 - u_E - u) (u_E + v_E + u + v)-Q (v_E + v)\big), \\
\dfrac{dv}{d\tau} &=(v_E + v) \big(C (u_E + u) - (u_E + v_E + u + v) (M + N (v_E + v))\big)\,.
\end{aligned}
\end{equation}
In particular, the Jacobian matrix of~\eqref{eq:Y} at the equilibrium point $(0,0)$ at the bifurcation point $(C^*,Q^*)$ coincides with $\partial {X}/\partial(u,v)\left(u_E,v_E;C^*,Q^*\right)$ in ~\eqref{eq:jac}.

Let $\mathbf{P}=[\mathbf{v}_1, \mathbf{v}_2]$ be the matrix whose columns are $\mathbf{v}_1$ and $\mathbf{v}_2$; see~\eqref{eq:eigenvectors}--(\ref{eq:v2}). 
Next, consider the following change of coordinates:
\begin{equation*}
 \left(
\begin{array}{c}
 x\\
 y \\
\end{array}
\right)
=
\mathbf{P}^{-1}\left(
\begin{array}{c}
 u\\
 v\\
\end{array}
\right).
\end{equation*}
Then, the vector field given by
\begin{equation*}
\mathbf{J}=\mathbf{P}^{-1}\circ Y\circ\mathbf{P},
\end{equation*}
is $\mathcal{C}^{\infty}$-conjugated to system~\eqref{eq:Y}.

Taking a Taylor expansion of $\mathbf{J}(x,y;C,Q)$ with respect to $(x,y)$ around $(x,y)=(0,0)$ and evaluating at $(C,Q)=(C^{\ast},Q^{\ast})$, one obtains
\begin{align*}\dfrac{d}{d\tau}
\begin{pmatrix}
x\\[1.6mm]y
\end{pmatrix}
=
\begin{pmatrix}
0 & 1 \\[1.6mm]  0 & 0
\end{pmatrix}
\begin{pmatrix}
x \\[1.6mm] y
\end{pmatrix}
+
\dfrac{1}{z_1 z_2-1}
\begin{pmatrix}
a_{20}x^2+a_{11}xy+a_{01}y^2+\mathcal{O}(||(x,y)||^3) \\[1.6mm]
b_{20}x^2+b_{11}xy+b_{02}y^2+\mathcal{O}(||(x,y)||^3)
\end{pmatrix}
\end{align*}
provided $z_1 z_2-1\neq0$,
where 
\begin{align*}
a_{20}=&2( C^* z_1 -M - N u_E - M z_1 - 3 N v_E - 2 N z_1 v_E - z_1 z_2+Q^* z_1 z_2+ 2 u_E z_1 z_2 - z_1^2 z_2 + 3 v_E z_1^2 z_2\\
& + z_1^2 v_E z_2);\\ 
b_{20}	=&-2z_1G_2;~\text{and}\\ 
b_{11}	=&1 - Q^* - 2 u_E + 2 z_1 - C^* z_1 + M z_1 - 6 u_E z_1 - 2 z_1 v_E+2 N z_1 v_E + z_1 z_2+ 2 M z_1 z_2 - Q^* z_1 z_2 \\ 
 		&- 2 u_E z_1 z_2 +2 N u_E z_1 z_2 - C^* z_1^2 z_2 + M z_1^2 z2+ 6 N v_E z_1 z_2 + 2 N z_1^2 v_E z_2.
\end{align*}

If $b_{20}\neq0$ and $a_{20}+b_{11}\neq0$, then the theory of normal forms for bifurcations~\cite{kuznetsov} ensures that our system fulfils the necessary genericity conditions to undergo a codimension two Bogdanov-Takens bifurcation. 
In particular, since $z_1\neq0$ in~\eqref{eq:v1}, condition
\begin{equation}\label{eq:bt2}
G_2\neq0
\end{equation}
ensures that $b_{20}\neq0$. Furthermore, after some algebraic manipulation one obtains
$ 
 a_{20}+b_{11}= G_3G_4.
$
Hence, condition $a_{20}+b_{11}\neq0$ is equivalent to 
\begin{equation}\label{eq:bt3}
G_3G_4\neq0.
\end{equation}

In summary, step 1 and inequalities (\ref{eq:bt1}), (\ref{eq:bt2}) and (\ref{eq:bt3}) ensure that the genericity and transversality conditions of a codimension two Bogdanov-Takens normal form are satisfied. Hence, there exists a smooth, invertible transformation of coordinates, an orientation-preserving time rescaling, and a reparametrisation such that, in a sufficiently small neighbourhood of $(u,v,C,Q)=(u_E,v_E,C^{\ast},Q^{\ast})$, the system (\ref{model2}) is topologically equivalent to one of the following normal forms of a Bogdanov-Takens bifurcation:
\begin{equation}
\begin{array}{rcl} \label{eq:nf-bt}
  \dot{\xi}_1 & = &\xi_2, \\
  \dot{\xi}_2 & = & \beta_1+\beta_2 \xi_2 +\xi_2^2\pm \xi_1\xi_2,\\
\end{array}
\end{equation}
where the sign of the term $\xi_1\xi_2$ in (\ref{eq:nf-bt}) is determined by the sign of $b_{20}(a_{20}+b_{11})$.

We are now in a position to state the corresponding result.

\begin{theorem}\label{btp3}
Let the system parameter be such that $\Delta=0$~\eqref{delta}, $tr(J(E))=0$~\eqref{qeq} and $G_{1,2,3,4}\neq0$~\eqref{eq:G1},~\eqref{eq:G2} -- \eqref{eq:G4}, then the equilibrium point $E=(u_E,v_E)$ in (\ref{model2}) undergoes a codimension-two Bogdanov-Takens bifurcation.		
\end{theorem}

\subsection{Bifurcation diagram}
\begin{figure}
\centering
\includegraphics[width=9cm]{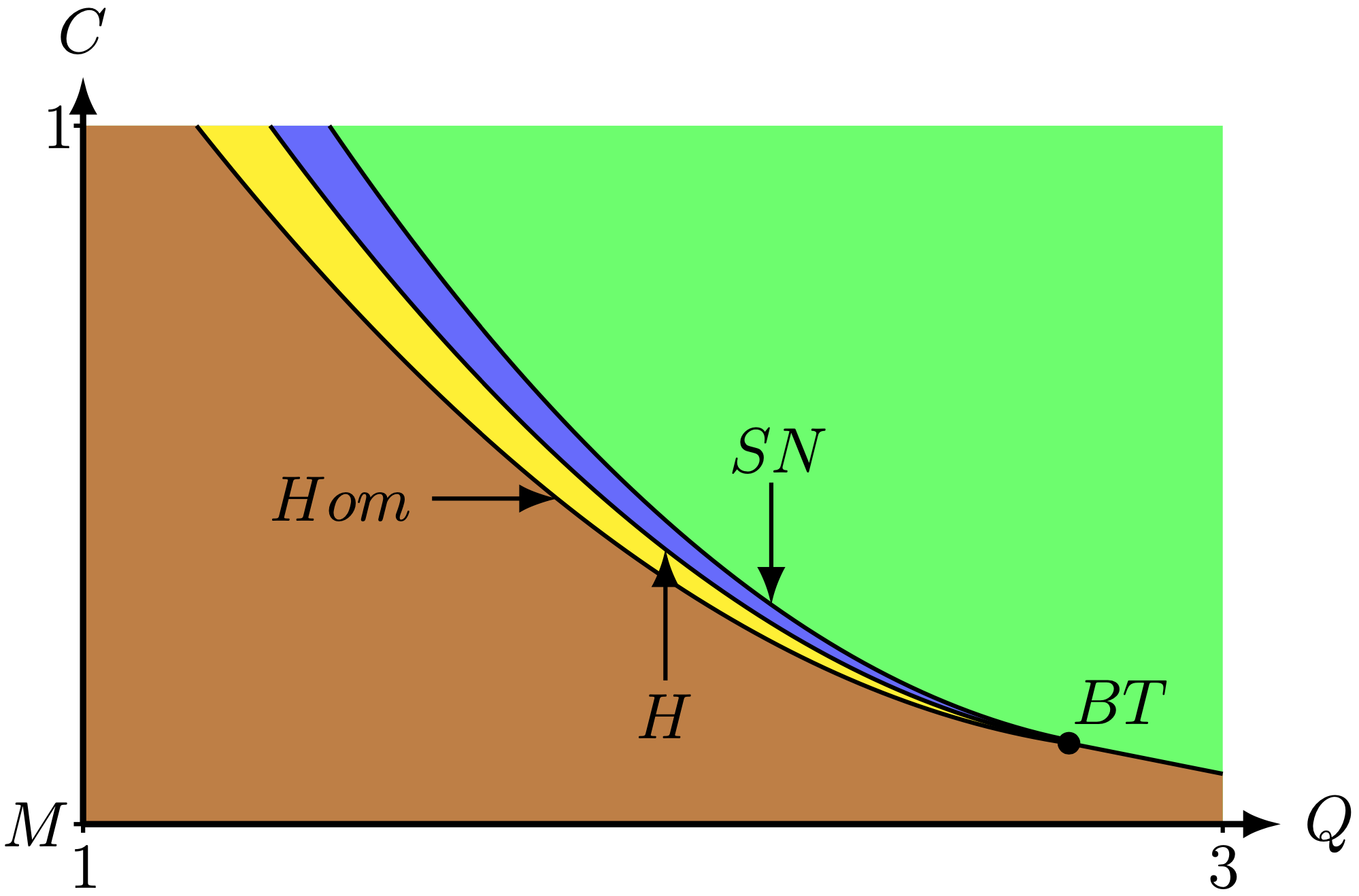}
\caption{The bifurcation diagram of system~\eqref{model2} for $(M,N)=(0.16,0.25)$ fixed and created with the numerical bifurcation package MATCONT. The curve $H$ represents the Hopf curve~\eqref{sec:hopf}, $SN$ represents the saddle-node curve~\eqref{sn}, $Hom$ represents the homoclinic bifurcation and $BT$ represents the Bogdanov–Takens bifurcation~\eqref{sec:bt}. See Figure~\ref{Fig10} for typical phase planes in the four different regions. }
\label{Fig9}
\end{figure}

In order to obtain the bifurcation diagram of system~\eqref{model2} for $Q>1$, $C>M$ and the parameters $M$ and $N$ fixed we follow~\cite{arancibia4,arancibia5} and use the numerical bifurcation package MATCONT~\cite{matcont}\footnote{Note that the Matlab package ode45 was used to generate the data for the simulations and then the PGF package (or tikz) was used to generate the graphics format.}. The bifurcation curves obtained from Sections ~\ref{SN_S} -- \ref{sec:bt} divide the $(Q,C)$-parameter space into four parts, see Figure~\ref{Fig9}. 
We observe that system~\eqref{model2} has only one positive equilibrium point if $(Q,C)$ are located on the saddle-node curve, see Theorem~\ref{p3}, 
while system~\eqref{model2} does not have any equilibrium points in the first quadrant, and therefore $(0,0)$ is global attractor, if $(Q,C)$ are located above the saddle-node curve (green region of Figure~\ref{Fig9}). 
System~\eqref{model2} has two positive equilibrium points, namely $P_1=\left(u_1,v_1\right)$ and $P_2=\left(u_2,v_2\right)$ with $0<u_1\leq u_2$ and $0<v_1\leq v_2$, if $(Q,C)$ are located below the saddle-node curve (blue, yellow and brown regions of Figure~\ref{Fig9}). In these regions the equilibrium point $P_1$ is always a saddle point. 
The equilibrium point $P_2$ is stable when $(Q,C)$ are located below the Hopf curve (yellow and brown regions of Figure~\ref{Fig9}) and it is surrounded by an unstable limit cycle when $(Q,C)$ is, in addition, above the homoclinic curve (yellow region of Figure~\ref{Fig9}).
The equilibrium point $P_2$ is unstable when $(Q,C)$ are located between the Hopf curve and the saddle-node curve (blue region Figure~\ref{Fig9}). See Figure~\ref{Fig10} for typical phase planes each of the four regions.

\section{Conclusions}\label{con}

In this manuscript, a Bazykin predator-prey model with predator intraspecific interactions and the ratio-dependent functional response was studied. Using a diffeomorphism we transformed the Bazykin predator-prey model~\eqref{model1} to a topologically equivalent system~\eqref{model2} and subsequently analysed this nondimensionalised system. 
\begin{figure}
\centering
\includegraphics[width=12.5cm]{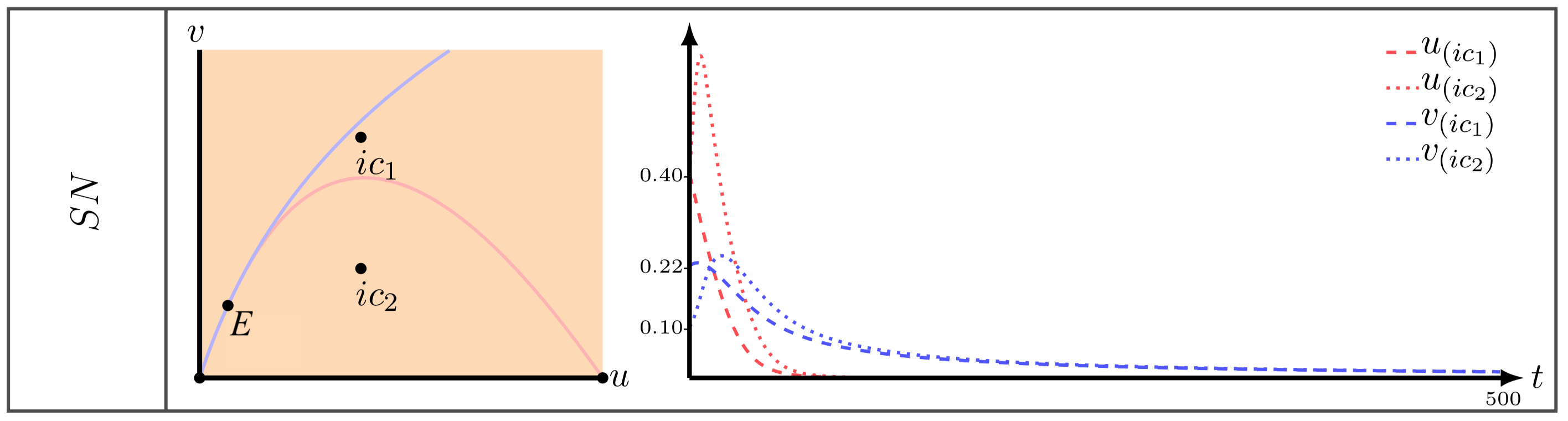}
\includegraphics[width=12.5cm]{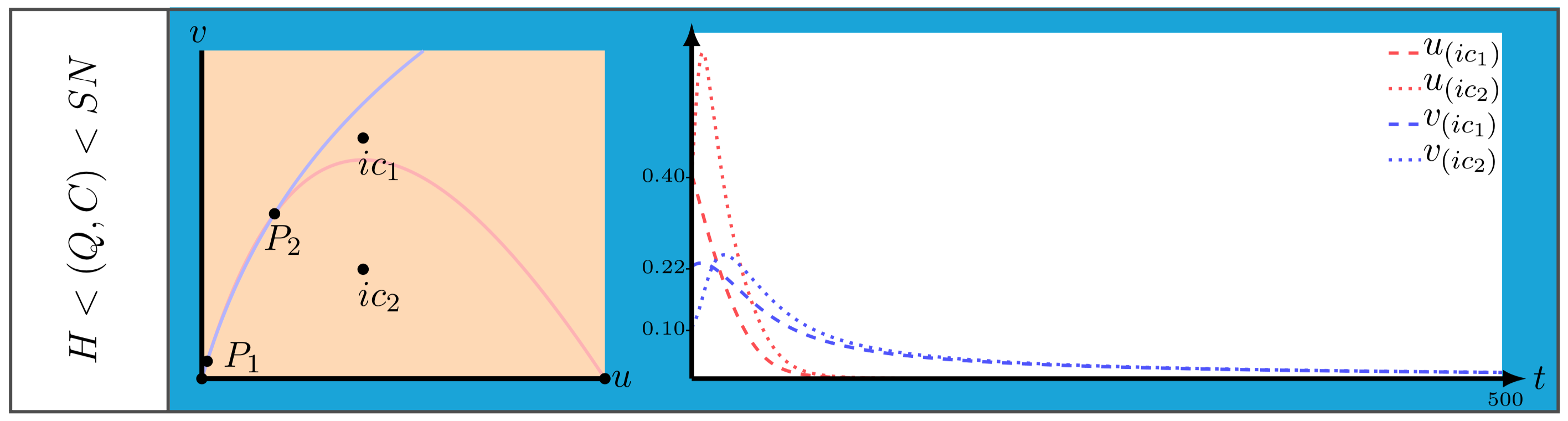}
\includegraphics[width=12.5cm]{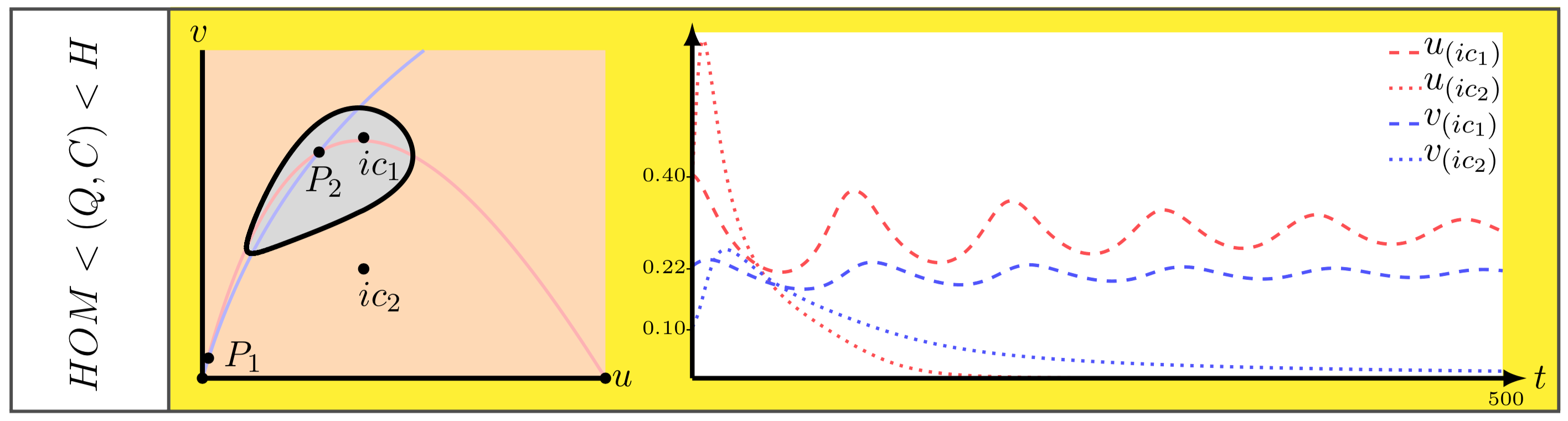}
\includegraphics[width=12.5cm]{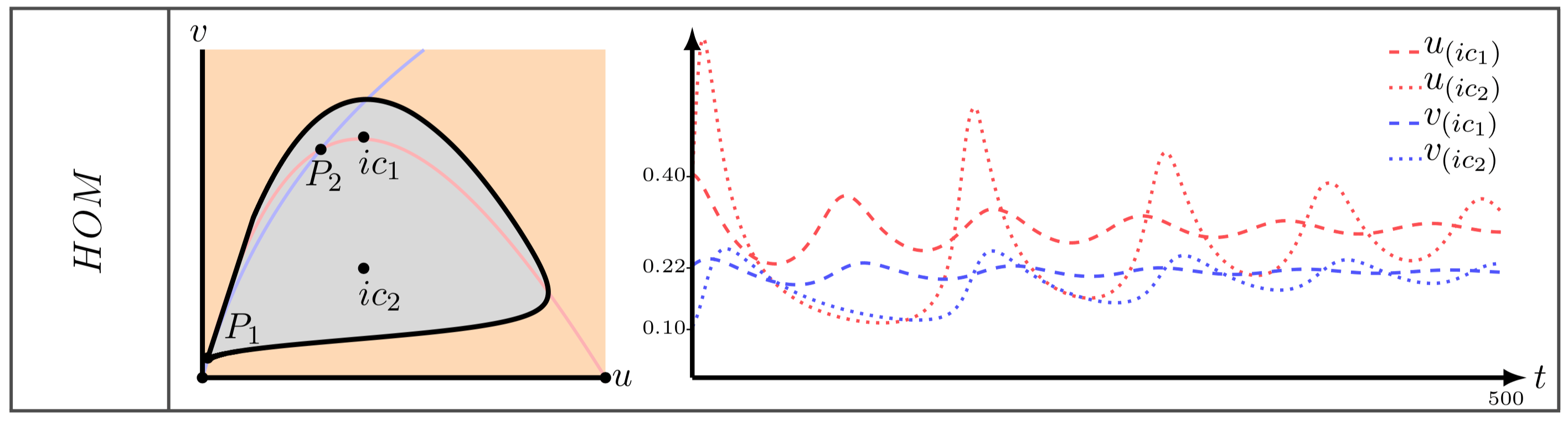}
\includegraphics[width=12.5cm]{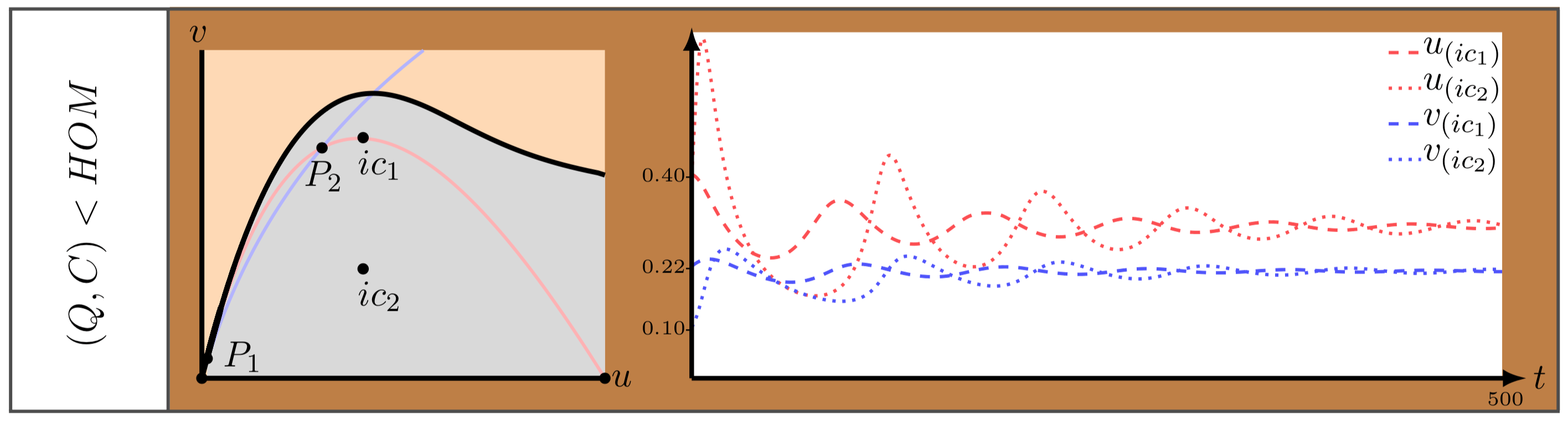}
\caption{The characteristics of predator-prey populations in the ratio-dependent Bazykin model~\eqref{model2} over time ($t$) for the $(Q,C)$-parameter space showed in Figure~\ref{Fig9}. 
In the left panels we show the phase planes for the system parameters $(C,M,N)=(0.363,0.16,0.25)$ fixed and decreasing $Q$, with $Q=1.83$ in the top  panel where a saddle-node bifurcation occurs. In the right panels we show the time series dynamics obtained through numerical integration considering the initial densities $ic_{1}=\left(u_{\left(ic_1\right)},v_{\left(ic_1\right)}\right)=\left(0.4,0.22\right)$ and $ic_{2}=\left(u_{\left(ic_2\right)},v_{\left(ic_2\right)}\right)=\left(0.4,0.1\right)$.}
\label{Fig10}
\end{figure}
This system has four system parameters which determine the number of the equilibrium points and their stability. We showed that the equilibrium point $P_1$ is always a saddle point and $(1,0)$, which corresponds to the rescaled carrying capacity, is a saddle point if $C>M$ and a stable node if $C<M$. We showed in Theorem~\ref{00} that the origin has complex dynamics and by using vertical and horizontal blow-ups we determined the dynamic in the neighbourhood of the origin. Furthermore, for some sets of parameters values the stable manifold of $P_1$ determines a separatrix curve which divides the basins of attraction of $(0,0)$ and $P_2$. As a result, the equilibrium point $P_2$ can be stable, stable surrounded by unstable limit cycle or unstable, depending on the trace of its Jacobian matrix, see Theorem~\ref{p2}. Additionally, we proved that system~\eqref{model2} undergoes a codimension-one Hopf bifurcation at the equilibrium $P_2$ and we also provided evidence for a degenerate codimension-two Hopf bifurcation at the Bautin point, see Theorem~\ref{hopfp3}. Moreover, the equilibrium points $P_1$ and $P_2$ collapse for $\Delta=0$~\eqref{delta}, i.e. $P_1=P_2=E$, and system~\eqref{model2} undergoes to a saddle-node bifurcation~\cite{perko}, see Theorem~\ref{sn}, or a Bogdanov-Takens bifurcation, see Theorem~\ref{btp3}.

Since the function $\varphi$~\eqref{diff} is a diffeomorphism preserving the orientation of time, the dynamics of system~\eqref{model2} is topologically equivalent to the dynamics of Bazykin predator-prey model~\eqref{model1}. Therefore, we can conclude that there exists self-regulation in system~\eqref{model1} for certain population sizes, that is, the species can coexist. However, system~\eqref{model2} is sensitive to changes in the parameters and also disturbances of the population size. We can see this impact in the size of the basins of attraction of the equilibrium points $(0,0)$ and $P_2$ in Figures~\ref{Fig5} and~\ref{Fig10}. In these figures, the orange region represents population sizes that lead to the extinction of both populations and the grey region represents population sizes that lead to the stabilisation of both populations over time. We showed that the stabilisation of the predator and the prey populations depend on the values of the parameters $Q$ and $C$ by taking the parameters $M$ and $N$ fixed\footnote{Note that the parameter $Q$ correspond to the rescaled per capita predation rate $q$ and the parameter $C$ correspond to the rescaled efficiency with which predators convert consumed prey into new predators $c$.}. 
For example, small predation rate $q$, or small efficiency with which predators convert consumed prey into new predators $c$, increases the area of coexistence which is related to basins of attraction of $P_2$, see Figure~\ref{Fig10}. Moreover, we show in Figure~\ref{Fig10} that the stabilisation or extinction of both the predator and the prey populations depends also on the species initial density. For example, for low predation rate the stable manifold of the equilibrium point $P_1$ act as a separatrix and initial conditions with large predator density lead to the extinction of both populations, see bottom panel of Figure~\ref{Fig10}. In contrast, if the parameters $(Q,C)$ are located in the yellow region of the bifurcation diagram, see Figures~\ref{Fig9} and~\ref{Fig10}, the stabilisation of both populations are bounded by the unstable limit cycle. In this case, the stabilisation or extinction depends on the initial conditions as is showed in Figure~\ref{Fig10}. If the initial condition is in the region interior to the limit cycle then the predator and prey population will coexist. However, if the initial condition is in the region outside of the limit cycle then both populations become extinct. 

Finally, we observe that intraspecific interaction occurs in the predator population.  Also, in most cases considered in this manuscript either both species go extinct or a coexistence state will be reached. For instance, in the left panel of Figure~\ref{Fig2} we observe the extinction of both populationswe observe the extinction of both populations (while in the right panel only the predator population goes extinct).  The dynamics observed in the left panel can be connected to the intraspecific interaction in the predator population. As observed in most of our results presented, this type of interaction improves the species’ adaption (see Figure~\ref{Fig5}). For example, when the prey is abundant, which may limit the predator population density, the intraspecies competition among predators becomes stronger. Hence, the predator density is insufficient to control the prey population and thus the prey is more likely to escape the predator~\cite{georgescu,bazykin}. 

\section*{Acknowledgment}
Pablo Aguirre was funded by Proyecto Interno UTFSM PI-LI-19-06. Pablo Aguirre also thanks Proyecto Basal CMM Universidad de Chile.

\newpage

\end{document}